\documentclass[12pt]{article}
\input{epsf.sty}
\usepackage{hyperref}
\usepackage{graphicx}
\usepackage{amsmath, amssymb}
\usepackage[arrow, matrix, curve]{xy}
\newtheorem {thm}{Theorem}[section]
\newtheorem {prop}[thm]{Proposition} 

\newtheorem {lem}[thm]{Lemma}
\newtheorem {cor}[thm]{Corollary}
\newtheorem {defn}[thm]{Definition}

\newtheorem {rem}[thm]{Remark}

\newcommand{\qed}{\nobreak \ifvmode \relax \else
      \ifdim\lastskip<1.5em \hskip-\lastskip
      \hskip1.5em plus0em minus0.5em \fi \nobreak
      \vrule height0.75em width0.5em depth0.25em\fi}

\def\Cox{\hfill \Box}

\def\N{{\Bbb N}}

\def\Z{{\Bbb Z}}
\def\R{{\Bbb R}}

\def\E{{\Bbb E}}

\def\D{\Delta}
\def\a{\alpha}

\def\b{\beta}

\def\d{\delta}

\def\phi{\varphi}

\def\l{\lambda}

\def\s{\sigma}

\def\x{\xi}

\def\o{\omega}

\def\D{\Delta}

\def\L{\Lambda}

\def\G{\Gamma}

\def\T{\T}

\def\FF{{\cal F}}

\def\AA{{\cal A}}

\def\PP{{\cal P}}

\begin{document}
\title{Gibbsian and non-Gibbsian properties of the generalized mean-field fuzzy Potts-model 
}

\author{
Benedikt Jahnel
\footnote{\scriptsize{Ruhr-Universit\"at   Bochum, Fakult\"at f\"ur Mathematik, D-44801 Bochum, Germany,
\newline
\texttt{Benedikt.Jahnel@ruhr-uni-bochum.de}, 
\newline
\texttt{www.ruhr-uni-bochum.de/ffm/Lehrstuehle/Kuelske/jahnel.html }}}
 \,, Christof K\"ulske
\footnote{\scriptsize{Ruhr-Universit\"at   Bochum, Fakult\"at f\"ur Mathematik, D-44801 Bochum, Germany,
\newline
\texttt{Christof.Kuelske@ruhr-uni-bochum.de}, 
\newline
\texttt{www.ruhr-uni-bochum.de/ffm/Lehrsttuehle/Kuelske/kuelske.html
/$\sim$kuelske/ }}}\, 
 \,, Elena Rudelli
\footnote{\scriptsize{Ruhr-Universit\"at   Bochum, Fakult\"at f\"ur Mathematik, D-44801 Bochum, Germany,
\newline
\texttt{Elena.Rudelli@ruhr-uni-bochum.de}}, 
}\, 
 \, \\
 and  Janine Wegener
\footnote{\scriptsize{Ruhr-Universit\"at   Bochum, Fakult\"at f\"ur Mathematik, D-44801 Bochum, Germany,
\newline
\texttt{Janine.Wegener@ruhr-uni-bochum.de}}, 
}\, 
\,  
}

\maketitle

\begin{abstract} We analyze the generalized mean-field $q$-state Potts model which is obtained 
by replacing the usual quadratic interaction function in the mean-field Hamiltonian by a higher power $z$.  
We first prove a generalization of the known limit result for the empirical magnetization vector 
of Ellis and Wang \cite{ElWa89} which shows that in the right parameter regime, the first-order phase-transition persists. 

Next we turn to the corresponding generalized fuzzy Potts model which is obtained 
by decomposing the set of the $q$ possible spin-values into $1<s<q$ classes and 
identifying the spins within these classes. 
In extension of earlier work \cite{HK04} which treats the quadratic model we prove the following: 
The fuzzy Potts model with interaction exponent bigger than four (respectively bigger than two and smaller or equal four) 
is non-Gibbs if and only if its inverse temperature $\b$ satisfies 
$\b\geq \b_c(r_*,z)$ where $\beta_c(r_*,z)$ is the critical inverse temperature 
of the corresponding Potts model and $r_*$ is the size of the smallest class which is 
greater or equal than two (respectively greater or equal than three.)
 
We also provide a dynamical interpretation considering sequences of fuzzy Potts models 
which are obtained by increasingly collapsing classes at finitely many times $t$ and discuss 
the possibility of a multiple in- and out of Gibbsianness, depending on the collapsing scheme. 
\end{abstract}

\smallskip
\noindent {\bf AMS 2000 subject classification:} 82B20, 82B26.

 \smallskip
\noindent {\bf Keywords:}  Potts model, Fuzzy Potts model, Ellis-Wang Theorem,  
Gibbsian measures, non-Gibbsian measures, mean-field measures. 


%
%
%
%
%
%
%
%
%
%
%
%

\section{Introduction}

Past years have seen a number of examples of measures which arise from local transforms 
of Gibbs measures which turned out to be non-Gibbs, for a general background see  
\cite{EFS, DEZ, KUL6}. 
Two particularly interesting types of transformations which were considered recently 
are time-evolutions \cite{ACD1, AFHR10}
and local coarse-grainings \cite{B10, KULOP08, JaKu12},  
both without geometry (mean-field) and with geometry. 
Very recently in \cite{FHM13} there is even been considered a system of Ising spins on a large discrete torus with a Kac-type interaction subject to an independent spin-flip dynamics, using large deviation techniques (usually applied in the mean-field setting) for the empirical density allowing for a spatial structure with geometry.

In the present paper we pick up a line of a 
mean-field analysis which was begun in \cite{HK04}. 
The extension to exponents $z\geq 2$ is natural since it amounts to considering energy given 
by the number of $z$-cliques of equal color in the case of integer $z$, see \ref{Random cluster representation and $z$-clique variables}.
In \cite{HK04} the mean-field Potts model was
considered under a local coarse-graining. Here the local spin-space $\{1, \dots, q\}$ 
is decomposed into $1<s<q$ classes of sizes $r_1,\dots, r_s$. This map, performed at each 
site simultaneously, defines a coarse-graining map $T:\{1,\dots, q\}^N \rightarrow  \{1,\dots, s\}^N$. The measures arising 
as images of the Potts mean-field measures for $N$ spins under $T$ 
constitute the so-called fuzzy-Potts model. It was shown that non-Gibbsian behavior occurs if the temperature 
of the Potts model is small enough and precise transition-values between Gibbsian and non-Gibbs images 
were given. We remark that the notion of a Gibbsian mean-field model is employed which 
considers as a defining property the existence and continuity of single-site probabilities. This notion 
is standard by now (see for example \cite{KUL1,KUL2, FHM13-,JaKu13}) and provides the natural counterpart of Gibbsianness for lattice systems for mean-field measures. 

Aim one of the paper is to generalize the mean-field Potts Hamiltonian, and analyse 
phase-transitions for the generalized mean-field Potts measures. Is there an analogue of 
the Ellis-Wang theorem \cite{ElWa89} and persistance of the first order phase-transition? We show 
that this is indeed the case for $q>2$. For $q=2$ there is a threshold for the exponent such that for $2\leq z\leq4$ there is a phase-transition of second order, for $z>4$ the phase-transition is of first order.  

Aim two of the paper is to look at the Gibbsian properties of the resulting fuzzy model, obtained 
by application of the same map $T$ to the generalized mean-field Potts measure. 
Do we obtain the same characterization as for the standard mean-field Potts model? 
The answer is yes, but with changes, which are inherited by the changed behavior of the Curie-Weiss 
model when the interaction exponent changes. 

The third aim is to reinterpret our results and introduce a dynamical point of view. 
In this view we consider decreasing finite sequences of decompositions $\AA_t$, of the local state-space $\{1,\dots,q\}$, 
labelled by a discrete time $t=0,1,\dots,T$.  We call these sequences collapsing schemes. 
As we move along $t$ 
we are interested in whole trajectories of fuzzy measures and what can be said about Gibbsianness here. 
Analogous questions have been studied for time-evolved Gibbs measures arising from stochastic spin dynamics and usually 
there is no multiple in- and out Gibbsianness in these models. 
As we will see this may very well be the case here, depending on the collapsing scheme. 

Technically the paper rests on a detailed bifurcation analysis of the free energy, the first 
step being a reduction to a one-dimensional problem using an extension of the 
proof of \cite{ElWa89}. We find here 
the somewhat surprising fact  that there is a triple point for  $q=2,z=4$, with a transition from second-order to first-order phase-transition.

%

\section{The generalized Potts model}\label{The generalized Potts model}

For a positive integer $q$ and a real number $z\geq2$, the Gibbs measure $\pi^N_{\b,q,z}$ for the $q$-state generalized Potts model on the complete graph with $N\in\N$ vertices at inverse temperature $\b\geq0$, is the probability measure on $\{1,\dots, q\}^N$ which to each $\xi\in\{1,\dots,q\}^N$ assigns probability 
\begin{equation}\label{Gen_Potts_Finite}
\begin{split}
\pi^N_{\b,q,z}(\xi)=\frac{1}{Z^N_{\b,q,z}}\exp(-NF_{\b,q,z}(L_N^\xi))
\end{split}
\end{equation}
where $L_N^\xi=\frac{1}{N}\sum_{i=1}^N1_{\xi_i}$ is the empirical distribution of the configuration $\xi=(\xi_i)_{i\in N}$, $F_{\b,q,z}:\PP(\{1,\dots,q\}^N)\to\R$, $F_{\b,q,z}(\nu):=-\b\sum_{i=1}^q\nu_i^z/z$ is the mean-field Hamiltonian of the generalized Potts model and $Z^N_{\b,q,z}$ is the normalizing constant. Notice, the case $z=2$ is the standard Potts model, in particular the case $z=2$, $q=2$ refers to the Curie-Weiss-Ising model. We call the case $q=2$, $z\geq2$ the generalized Curie-Weiss-Ising model.

The Ellis-Wang Theorem \cite{ElWa89} describes the limiting behaviour of the standard Potts model as the system size grows to infinity. Here we give a generalized version for interactions with $z\geq2$.
\begin{thm}(Generalized Ellis-Wang Theorem)\label{Generalized_Ellis_Wang}
Assume that $q\geq2$ and $z\geq2$ then there exists a critical temperature $\b_c(q,z)>0$ such that in the weak limit
\begin{equation}
\lim_{N\to\infty}\pi^N_{\b,q,z}(L_N\in\cdot)
=\begin{cases}
\d_{1/q( 1,\dots,1)} & \text{ if } \b<\b_c(q,z)\cr
\frac{1}{q}\sum^q_{i=1}\d_{u(\b,q,z)e_i+(1-u(\b,q,z))/q(1,\dots, 1)}& \text{ if } \b>\b_c(q,z)
\end{cases}
\end{equation}
where $e_i$ is the unit vector in the $i$-th coordinate of \textit{ }$\R^q$ and $u(\b,q,z)$ is the largest solution of the so-called mean-field equation
\begin{equation}\label{MFeqPhi}
u=\frac{1-\exp\bigl(\D_{\b,q,z}(u)\bigr)}{1+\left(q-1\right)\exp\bigl(\D_{\b,q,z}(u)\bigr)}
\end{equation}
with $\D_{\b,q,z}(u):=-\frac{\beta}{q^{z-1}}\left[\bigl(1+(q-1)u\bigr)^{z-1}-\bigl(1-u\bigr)^{z-1}\right]$. 

Further, for $(q,z)\in\{2\}\times[2,4]$ the function $\b\mapsto u(\b,q,z)$ is continuous. In the complementary case the function $\b\mapsto u(\b,q,z)$ is discontinuous at $\b_c(q,z)$.
\end{thm}

For $q>2$ the above result is in complete analogy to the standard Potts model. For the generalized Curie-Weiss-Ising model ($q=2$) there is an important difference.
It is a known fact that the standard Curie-Weiss-Ising model ($z=2$) has a second order phase-transition. This is still true as long as $2\leq z\leq4$.
But in case of the generalized Curie-Weiss-Ising model with $z>4$ the phase-transition is of first order.


\bigskip
In the analysis of the fuzzy Potts model the following result is useful.
\begin{prop}\label{q_Monotonicity_Of_Beta}
For the generalized Potts model the function $q\mapsto\b_c(q,z)$ is increasing.
\end{prop}

\section{The generalized fuzzy Potts model}\label{The generalized fuzzy Potts model}
Consider the $q$-state generalized Potts model and let $s<q$ and $r_1,..., r_s$ be positive integers such that $\sum^s_{i=1}r_i=q$. For fixed $\b>0$, $z\geq2$ and $N\in\N$ let $X$ be the $\{1,\dots, q\}^N$-valued random vector distributed according to the Gibbs measure $\pi^N_{\b,q,z}$. Then define $Y$ as the $\{1,\dots, s\}^N$-valued random vector by
\begin{align*}
Y_i=
	\begin{cases}
	1\quad \text{if}\quad X_i\in\{1,..., r_1\},\\
	2\quad \text{if}\quad X_i\in\{r_1+1,..., r_1+r_2\},\\	
	\vdots \quad\quad\quad \vdots\\
	s\quad \text{if}\quad X_i\in\{q-r_s+1,..., q\}
	\end{cases}
\end{align*}
for each $i\in\{1,\dots,N\}$. In other words using the coarse-graining map $T:\{1,\dots,q\}^N\mapsto\{1,\dots,s\}^N$ with $T(k)=l$ iff $\sum_{j=1}^{l_i-1}r_{j}<k_i\leq\sum_{j=1}^{l_i}r_{j}$ for all $i\in\{1,\dots,N\}$ we have $Y=T\circ X$. Let us denote $\mu^N_{\b,q,z,(r_1,...,r_s)}$ the distribution of $Y$
and call it the finite-volume fuzzy Potts measure. The vector $(r_1,..., r_s)$ we call the spin partition of the fuzzy Potts model.

\bigskip
In \cite{HK04} the notion of Gibbsianness for mean-field models is introduced. It is based on the continuity of the so-called mean-field specification as a function of the boundary condition. In analogy to the lattice situation a mean-field specification is a probability kernel that for every boundary measure is a measure on the single site space. If it is discontinuous w.r.t the boundary measure, it cannot constitute a Gibbs measure. The mean-field specification is obtained as the infinite-volume limit of the one site conditional probabilities in finite volume. To be more specific we present the statement from \cite{HK04} applied to our situation without proof.
\begin{lem}
For $\mu^N_{\b,q,z,(r_1,...,r_s)}$ the generalized fuzzy Potts model on $\{1,\dots,s\}$ there exists a probability kernel $Q^N_{\b,q,z,(r_1,...,r_s)}: \{1,\dots,s\}\times \PP(\{1,\dots,s\})\to[0,1]$ such that the single-site conditional expectations at any site 
$i$ can be written in the form 
\begin{equation*}
\mu^N_{\b,q,z,(r_1,...,r_s)}(Y_i=k|Y_{\{1,\dots,N\}\setminus i}=\eta)=Q^N_{\b,q,z,(r_1,...,r_s)}(x|\bar\eta)
\end{equation*}
where $\bar\eta\in\PP(\{1,\dots,s\})$ with $\bar\eta_l=\#(1\leq j\leq N,j\neq i,\eta_j=l)/(N-1)$ the fraction of sites for which the 
spin-values of the conditioning are in the state $l\in\{1,\dots,s\}$. 
Further $\mu^N_{\b,q,z,(r_1,...,r_s)}$ is uniquely determined by $Q^N_{\b,q,z,(r_1,...,r_s)}$.
\end{lem}
\begin{defn}
Assume for all $k\in\{1,\dots,s\}$ and $\nu_N\to\nu$, the infinite-volume limit $Q^N_{\b,q,z,(r_1,...,r_s)}(k|\nu_N)\to Q^\infty_{\b,q,z,(r_1,...,r_s)}(k|\nu)$ exists. We call the generalized fuzzy Potts model Gibbs if $\nu\mapsto Q^\infty_{\b,q,z,(r_1,...,r_s)}(\cdot|\nu)$ is continuous. Otherwise we call it non-Gibbs.
\end{defn}
Theorem 1.2 in \cite{HK04} therefore describes properties of the limiting conditional probabilities in case of the fuzzy Potts model. Here we give a version of this theorem for the generalized fuzzy Potts model with exponent $z>2$. 
\begin{thm}\label{Generalized_HaKu}
Consider the $q$-state generalized fuzzy Potts model at inverse temperature $\b>0$ with exponent $z>2$ and spin partition $(r_1,..., r_s)$, where $1<s<q$ and $\sum_{i=1}^s{r_i=q}$. 
Denote by $\beta_c(r_k,z)$ the inverse critical temperature of the $r_k$-state generalized Potts model with the same exponent $z>2$. Then 
\begin{itemize}
\item[(i)] Suppose $2<z\leq4$ and $r_i\leq2$ for all $i\in\{1,\dots,s\}$, then the limiting conditional probabilities exist and are continuous functions of empirical distribution of the conditioning for all $\b\geq0$.
\end{itemize}
Assume $z>4$ or that $r_i\geq3$ for some $i\in\{1,\dots,s\}$. Put $r_{*}:=\min\{r\geq 3,r=r_i \text{\ for some\ } i\in\{1,\dots, s\}\}$ and $r_{\#}:=\min\{r\geq 2,r=r_i \text{\ for some\ } i\in\{1,\dots, s\}\}$, then the following holds:
\begin{itemize}
\item[(ii)] If $z>4$ then
\begin{enumerate}
\item the limiting conditional probabilities exist and are continuous for all $\b<\b_c(r_{\#},z)$,
\item the limiting conditional probabilities are discontinuous for all $\b\geq\b_c(r_{\#},z)$, in particular they do not exist in points of discontinuity.
\end{enumerate}
\item[(iii)] If $2<z\leq 4$ then
\begin{enumerate}
\item the limiting conditional probabilities exist and are continuous for all $\b<\b_c(r_{*},z)$,
\item the limiting conditional probabilities are discontinuous for all $\b\geq\b_c(r_{*},z)$, in particular they do not exist in points of discontinuity.
\end{enumerate}
\end{itemize}
\end{thm}

\section{Dynamical Gibbs-non Gibbs transitions along collapsing schemes}

Consider the set of Potts spin values $\{1,\dots,q\}$ and denote by $\AA=\{I_1,\dots,I_r\}$ a spin partition.  
Write $\mu_{\b,q,z,\AA}^N$ for the finite-volume fuzzy Potts Gibbs measure 
on $\{1,\dots,r \}^N$. With a partition $\AA$ comes the $\s$-algebra 
$\s(\AA)$ which is generated by it. 
 It consists of the unions of sets 
in $\AA$. 
Conversely a $\s$-algebra determines 
a partition. 

The set of $\s$-algebras over $\{1,\dots,q\}$ is partially ordered by inclusion. 
Now let $(\AA_t)_{t=0,1,\dots, T}$ be a strictly decreasing sequence of partitions (a \textit{collapsing scheme}) with $\AA_0=(\{1\},\dots, \{q\})$ being the finest 
one (consisting of $q$ classes), and $\AA_T=(\{1,\dots,q\} )$ being the coarsest one.  
$t$ can be considered as a time index. Moving along $t$ more and more classes are collapsed. 
Note that the finite sequence of $\s$-algebras generated by these 
partitions, $\s(\AA_T)\subset \s(\AA_{T-1}) \subset \dots \subset \s(\AA_0)$ is a filtration. 
Such a filtration can be depicted as a rooted tree with $q$ leaves which has  
$T$ levels. A level $i$ corresponds to a $\s$-algebra $\FF_i$, 
the vertices at level $i$ are the sets in the partition corresponding to $\FF_i$. 
A set in the partition at level $i$ is a parent of a set in the partition at level $i-1$ iff 
it contains the latter.

We look at the corresponding sequence of increasingly coarse-grained models 
$(\mu_{\b,q,z,\AA}^N)_{t=0,\dots, T}$. 
What can be said about in and out of Gibbsiannes along such a path?
For a partition $\AA$ and given exponent $z\geq2$ denote by $r_*(\AA,z)$ the size of the smallest class in the non-Gibbsian region $(r,z)\in([2,\infty)\times[2,\infty))\setminus(\{2\}\times[2,4])$. The following corollary is a direct consequence of our main Theorem \ref{Generalized_HaKu} and Theorem 1.2 in \cite{HK04}.
\begin{cor}\label{Non_Gibbs_Cor}The model is non-Gibbs at time $t\in \{1,2,\dots, T-1\}$ if and only if $\b\geq\b_c(r_*(\AA_t,z),z)$. 
\end{cor}
Even though by Proposition \ref{q_Monotonicity_Of_Beta} $r\mapsto\b_c(r,z)$ is increasing, it is quite possible to have collapsing schemes where $t\mapsto\b_c(r_*(\AA_t,z),z)$
is not monotone for $t\in \{1,\dots, T\}$.
This is because $t\mapsto r_*(\AA_t,z)$ does not have to have monotonicity, as it happens e.g in the following example: 
\begin{equation*}\label{Rotation_Finite_Volume}
\begin{split}
&\AA_0=(\{1\},\{2\},\{3\},\{4\},\{5\})\cr
&\AA_1=(\{1,2\},\{3\},\{4\},\{5\})\cr
&\AA_2=(\{1,2,3\},\{4\},\{5\})\cr
&\AA_3=(\{1,2,3\},\{4,5\})\cr
&\AA_4=(\{1,2,3,4,5\})\cr
\end{split}
\end{equation*}
with $(r_*(\AA_t,5))_{t=1,\dots, T-1}=(2,3,2)$. 
If $q$ is a power of two, and the collapsing scheme is chosen according to a binary tree,  
there is of course monotonicity, as e.g. in the following example
\begin{equation*}\label{Rotation_Finite_Volume}
\begin{split}
&\AA_0=(\{1\},\{2\},\{3\},\{4\},\{5\},\{6\},\{7\},\{8\})\cr
&\AA_1=(\{1,2\},\{3,4\},\{5,6\},\{7,8\})\cr
&\AA_2=(\{1,2,3,4\},\{5,6,7,8\})\cr
&\AA_3=(\{1,2,3,4,5,6,7,8\})\cr
\end{split}
\end{equation*}
with 
$(r_*(\AA_t,5))_{t=1,\dots, T-1}=(2,4)$.
\begin{defn} Let us agree to call a collapsing scheme regular 
if and only if $(r_*(\AA_t,z))_{t=1,\dots, T-1}$ is increasing, $T\geq 2$ (meaning there is no immediate collapse.)  
\end{defn}
The following theorem is an immediate consequence of Corollary \ref{Non_Gibbs_Cor} and Proposition \ref{q_Monotonicity_Of_Beta}.
\begin{thm} Consider the generalized $q$-state Potts model 
with interaction exponent bigger than $2$. 
For a regular collapsing scheme the following is true: 
\begin{itemize}
\item[(i)] The model stays Gibbs forever iff $\b <\b_c(r_*(\AA_1,z),z)$. 
\item[(ii)] It is non-Gibbs for all $t\in \{1,\dots, T-1\}$ iff $\b \geq  \b_c(r_*(\AA_{T-1},z),z)$. 
\item[(iii)] For $\b\in (\b_c(r_*(\AA_{1},z),z),\b_c(r_*(\AA_{T-1},z),z)]$ there is a transition time $t_G\in \{2,\dots, T-1\}$ such that 
the model is non-Gibbs for $t \in \{1,\dots, t_G -1\}$ and Gibbs for $t \in \{t_G, \dots, T\}$. 
\end{itemize}
\end{thm}
Note that the second temperature-regime of non-Gibbsianness contains temperatures 
which are strictly bigger than the critical temperature of the initial $q$-state Potts model.  
In the last regime there is an immediate out of Gibbsiannes, then the model stays non-Gibbs for a while and becomes Gibbsian again at the transition time $t_G$.  
Also note that for general collapsing schemes there can be temperature regions for which 
multiple in and out of Gibbsianness will occur. 
\begin{figure}[h]
\begin{center}
\includegraphics[width=14cm]{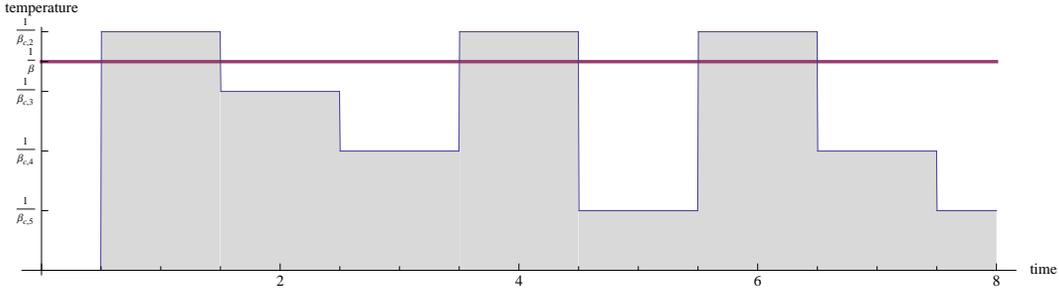}
\end{center}
\caption{\scriptsize{Qualitative picture of the collapsing scheme $(r_*(\AA_t,5))_{t=1,\dots, 7}=(2,3,4,2,5,2,3,5)$. The gray area below the graph shows the non-Gibbsian temperature regime. Clearly the  generalized Potts model with fixed temperature $1/\b$ and the same exponent under fuzzyfication given by the collapsing scheme $\AA_t$ can experience in and out of Gibbsianness multiple times.}}
\label{Collapsing scheme}
\end{figure}



\section{Proofs of statements presented in Section \ref{The generalized Potts model}}

\subsection{Proof of Theorem \ref{Generalized_Ellis_Wang}}
The empirical distribution $L_N$ obeys a large deviation principle with the relative entropy $I(\cdot|\a)$ as a rate function, where $\a$ is the equidistribution on $\{1,\dots,q\}$. Together with Varadhan's lemma the question of finding the limiting distribution of $L_N$ under $\pi_{\b,q,z}^N$ is equivalent to finding the global minimizers of the so-called free energy $\G_{\b,q,z}: \PP(\{1,\dots,q\})\mapsto\R$, 
\begin{equation}
\begin{split}
\G_{\b,q,z}(\nu)&=F_{\b,q,z}(\nu)+I(\nu|\a)=\frac{\b}{z}\sum_{i=1}^q\nu_i^z+\sum_{i=1}^q\nu_i\log(q\nu_i).
\end{split}
\end{equation}
For details on large deviation theory check \cite{DeZe10}. The proof of \ref{Generalized_Ellis_Wang} thus rests completely on the following theorem.
\begin{thm}\label{Minimizer_Theorem}
\begin{enumerate}
\item[(i)] Any global minimizer of $\G_{\b,q,z}$ with $z\geq2$ and $q\geq2$ must have the form
\begin{equation}
\bar{\nu}=
\begin{pmatrix}
\frac{1}{q}\bigl(1+\left(q-1\right)u\bigr)\\
\frac{1}{q}\left(1-u\right)\\
\vdots\\
\frac{1}{q}\left(1-u\right)
\end{pmatrix}
\quad \textup{\textit{with}} \quad u \in \left[0,1\right)
\end{equation}
or a point obtained from such a $\bar{\nu}$ by permutating the coordinates.
\item[(ii)] There exists a critical temperature $\beta_c(q,z)>0$ such that for $\beta<\beta_c(q,z)$, $\bar{\nu}$ is given in the above form with $u=0$, in other words $\bar{\nu}^T=(\frac{1}{q},\dots,\frac{1}{q})$. If $\beta>\beta_c(q,z)$, then $\bar{\nu}$ is given in the above form with $u=u(\b,q,z)$, where $u(\b,q,z)$ is the largest solution of the mean-field equation
\begin{equation}\label{MFeqPhi}
u=\frac{1-\exp\bigl(\D_{\b,q,z}(u)\bigr)}{1+\left(q-1\right)\exp\bigl(\D_{\b,q,z}(u)\bigr)}
\end{equation}
with $\D_{\b,q,z}(u):=-\frac{\beta}{q^{z-1}}\left[\bigl(1+(q-1)u\bigr)^{z-1}-\bigl(1-u\bigr)^{z-1}\right]$.
\item[(iii)] The function $\b\mapsto u(\b,q,z)$ is discontinuous at $\b_c(q,z)$ for all $z\geq2$ and $q\geq2$ except for the case $(q,z)\in\{2\}\times[2,4]$.
\end{enumerate}
\end{thm}

\bigskip
For the proof of part (i) of \ref{Minimizer_Theorem} we use the following remark and lemma. 
\begin{rem}\label{Permutation_Invariance}
Due to the permutation invariance of the model it suffices to consider minimizers of $\G_{\b,q,n}$ with $\bar{\nu}_{k}\geq\bar{\nu}_{k+1}$ for all $k\in\{1,\dots, q-1\}$.
\end{rem}
\begin{lem}\label{Only_Two_Poss}
Let $\bar{\nu}\in\PP(\{1,\dots,q\})$ be a minimizer of $\G_{\b,q,z}$ with $z\geq2$, $q\geq2$ and define the auxiliary function 
$$g(x):=\b x^{z-1}-\log(qx)$$
with $x\in(0,1]$. Let $\tilde{u}$ be the minimizer of $g$, given by
\begin{equation}\label{MinG}
\tilde{u}:=\frac{1}{\sqrt[z-1]{\b(z-1)}},
\end{equation}
then the coordinates of $\bar{\nu}$ satisfy the following conditions:
\begin{enumerate}
\item[(i)] If $\bar{\nu}_1\leq\tilde{u}$, then $\bar{\nu}_k=\bar{\nu}_1$ for all $k \in\{2,\dots,q\}$ and any minimizer of $\G_{\b,q,z}$ has the form
\begin{equation*}
\bar{\nu}=(\frac{1}{q},\dots,\frac{1}{q})^T.
\end{equation*}
\item[(ii)] If $\bar{\nu}_1 > \tilde{u}$, then $\bar{\nu}_k \in\{\bar{\nu_0},\bar{\nu}_1\}$ for all $k \in\{2,\dots,q\}$ with $\bar{\nu}_1>\bar{\nu}_0$ and $g(\bar{\nu}_0)=g(\bar{\nu_1})$. In this case any minimizer of $\G_{\b,q,z}$ has the form
\begin{equation*}
\bar{\nu}=(\underbrace{\bar{\nu}_1,\dots,\bar{\nu}_1}_{{}l\text{ times}},\bar{\nu}_0,\dots,\bar{\nu}_0)^T \text{ with }\bar{\nu}_1=\frac{1-(q-l)\bar{\nu}_0}{l},
\end{equation*}
where $1\leq l\leq q$.
\end{enumerate}
\end{lem}

\textbf{Proof: }Since $\bar{\nu}$ is a minimizer $\nabla \G_{\b,q,z}(\bar{\nu})=(c,\dots,c)^T$. In other words $-\b \bar{\nu}_k^{z-1}+\log(q\bar{\nu}_k)+1=c$ for all $k\in\{1,\dots,q\}$ and hence 
\begin{equation*}
g(\bar{\nu}_{1})=\b\bar{\nu}_{1}^{z-1}-\log(q\bar{\nu}_1)=\b \bar{\nu}_{k}^{z-1}-\log(q\bar{\nu}_k)=g(\bar{\nu}_{k})
\end{equation*}
for all $k\in\{1,\dots,q\}$. The function $g$ has the following properties: $\lim_{x\to0}g(x)=+\infty$; $g(1)=\beta-\log(q)$; $g'(x)=\beta(z-1)x^{z-2)}-1/x$ and thus $g$ attains its unique extremal point in $\tilde x$; $g''(x)=\b(z-1)(z-2)x^{z-3}+\frac{1}{x^2}>0$ and hence $g$ is strictly convex with global minimum attained in $\tilde x$. 

As a consequence $g$ is injective on $(0,\tilde u]$ and hence if $\bar{\nu}_1\leq\tilde u$ by Remark \ref{Permutation_Invariance} $\bar{\nu}_k\leq\bar{\nu}_1$ and thus from $g(\bar{\nu}_1)=g(\bar{\nu}_k)$ for all $k$ it follows $\bar{\nu}_k=\bar{\nu}_1$ for all $k$. So $\bar{\nu}$ must be the equidistribution.

If $\bar{\nu}_1>\tilde u$, since $g$ is strictly convex, $\lim_{x\to0}g(x)=+\infty$ and $g(\bar{\nu}_1)=g(\bar{\nu}_k)$ we have $\bar{\nu}_k\in\{\bar{\nu}_0,\bar{\nu}_1\}$ for all $k$ where $\bar{\nu}_0<\bar{\nu}_1$ such that $g(\bar{\nu}_0)=g(\bar{\nu}_1)$. Consequently again by Remark \ref{Permutation_Invariance} $\bar\nu$ must have the following form 
\begin{equation}\label{Pre_Minimizers}
\bar{\nu}=(\underbrace{\bar{\nu}_1,\dots,\bar{\nu}_1}_{{}l\text{ times}},\bar{\nu}_0,\dots,\bar{\nu}_0)^T \text{ with }2\leq l \leq q.
\end{equation}
Since $\bar{\nu}$ is a probability measure $l\bar{\nu}_1+(q-l)\bar{\nu}_0=1$ and hence $\bar{\nu}_1=(1-(q-l)\bar{\nu}_0)/l$.
$\Cox$

\bigskip
\textbf{Proof of Theorem \ref{Minimizer_Theorem} part (i): }
First note, for $\bar{\nu}\in\PP(\{1,\dots,q\})$ a minimizer of $\G_{\b,q,z}$, $k\in\{2,\dots,q\}$ and
$$D_{\bar{\nu}}^{k}:=\{\nu \in \PP(\{1,\dots,q\}):\nu_x=\bar{\nu}_x\text{ for all }x\in\{2,\dots,q\} \setminus\{k\}\}$$
of course
$\min_{\nu\in D_{\bar{\nu}}^{k}}\G_{\b,q,z}(\nu)=\G_{\b,q,z}(\bar{\nu})$. Using this and the above Lemma \ref{Only_Two_Poss} for fixed $k$ we can set $a\in[0,1]$ such that $\sum_{i\neq 1,k}\bar{\nu}_i=1-a$ where $\bar{\nu}$ is a minimizer. Hence $\nu_1+\nu_k=a$ and for $\nu\in D^{k}_{\bar{\nu}}$, $\G_{\b,q,z}(\nu)$ has to be minimized as a function of the variable $\nu_1$ alone. We calculate
\begin{align*}
\frac{\partial \G_{\b,q,z}}{\partial \nu_1} 
&=-\b(\nu_1^{z-1}-(a-\nu_1)^{z-1})+\log\frac{\nu_1}{a-\nu_1}
\end{align*}
and thus have to analyse the inequality
\begin{equation*}
h_l(x):=\b(x^{z-1}-(a-x)^{z-1})\leq\log\frac{x}{a-x}=:h_r(x).
\end{equation*}
Notice $h_l$ and $h_r$ are both point symmetric at $x=a/2$ and $h_l(a/2)=0=h_r(a/2)$. In particular $a/2$ is a candidate for the Minimum of $\nu_1\mapsto\G_{\b,q,z}(\nu_1,\bar\nu_2,\dots,\bar\nu_{k-1},a-\nu_1,\bar\nu_{k+1},\dots,\bar\nu_q)$ and if it is $\nu_1=\nu_k$.
By point symmetry is suffices to look at $h_l$ and $h_r$ on the set $[a/2,a]$. Requiring $h'_l(a/2)=h'_r(a/2)$ is equivalent to 
\begin{equation*}
\frac{a}{2}=\frac{1}{\sqrt[z-1]{\b(z-1)}}=\tilde u.
\end{equation*}
Let us collect some further properties of $h_l$ and $h_r$: Both functions are convex on $[a/2,a)$; $\lim_{x\to a}h_r(x)=\infty$ and $h_l(a)=\b a^{z-1}<\infty$; $h''_l(a/2)=0=h''_r(a/2)$. Also 
\begin{equation*}
\begin{split}
h'''_l(a/2)=2\b(z-1)(z-2)(z-3)(a/2)^{z-4}\text{ and
} h'''_r(a/2)=4(a/2)^{-3}
\end{split}
\end{equation*}
so if $a/2=\tilde u$ some minor calculations show $h'''_l(a/2)=h'''_l(a/2)$ iff $z=4$. In particular for $z<4$, $h'''_l(a/2)<h'''_l(a/2)$ and for $z=4$ higher orders show the graph of $h_l$ close to $a/2$ is lower than the one of $h_r$. That is why we have to distinguish two cases with several subcases each.

\bigskip
First let $2\leq z\leq4$. We show, there is either one or no additional point $x\in(a/2,a]$ such that $h'_l(x)=h'_r(x)$. Let us write the temperature as a function of solutions of $h'_l(x)=h'_r(x)$,
\begin{equation}\label{Beta_For_h'}
\b_{z,a}(x)=\frac{a}{(z-1)((a-x)x^{z-1}+x(a-x)^{z-1})}.
\end{equation}
This function is strictly increasing, indeed $\b'_{z,a}>0$ is equivalent to 
\begin{equation}\label{Beta_Ableitung}
a(z-1)-zx-(xz-a)(\frac{a-x}{x})^{z-2}<0.
\end{equation}
Setting $y=(a-x)/x$ we can write this equivalently as
\begin{equation}\label{Beta_Ableitung_Replaced}
\begin{split}
a(z-1)-z\frac{a}{y+1}-(\frac{a}{y+1}z-a)y^{z-2}&<0\cr
(z-1)y-1-((z-1)-y)y^{z-2}&<0\cr
\end{split}
\end{equation}
where $x\mapsto y$, $(a/2,a]\mapsto [0,1)$ is bijective. Notice $z-1>y$ and $y^{z-2}\geq y^2$, hence 
\begin{equation*}
\begin{split}
(z-1)y&-1-((z-1)-y)y^{z-2}<(z-1)y-1-((z-1)-y)y^{2}\cr
&=y^3-(z-1)(y^2-y)-1<y^3-3(y^2-y)-1=(y-1)^3<0.\cr
\end{split}
\end{equation*}
But this is true and thus $\b'_{z,a}>0$ and for every $\b\leq\b_{z,a}(a/2)=\frac{1}{z-1}(\frac{a}{2})^{1-z}$ there is no $x\in(a/2,a]$ with $h'_l(x)=h'_r(x)$, for every $\b>\frac{1}{z-1}(\frac{a}{2})^{1-z}$ there is exactly one $x\in(a/2,a]$ with $h'_l(x)=h'_r(x)$.

\bigskip
Subcase one, let $a/2\leq\tilde u$. This is equivalent to $\b\leq\frac{1}{z-1}(\frac{a}{2})^{1-z}$ and hence $h'_r>h'_l$ on $[a/2,a)$, in particluar there can not be a $x\in[a/2,a)$ such that $h_l(x)=h_r(x)$ and $h_r>h_l$ on $[a/2,a)$. Due to point symmetry $a/2$ is the unique global minimum of the free energy as a function of the first variable $\nu_1$ on $D^k_{\bar\nu}$. In particular $\nu_1\leq\tilde u$ and thus by Lemma \ref{Only_Two_Poss} part (i), the free energy minimizer is the equidistribution.

\bigskip
Subcase two, let $a/2>\tilde u$. This is equivalent to $\b>\frac{1}{z-1}(\frac{a}{2})^{1-z}$ and hence there is exactly one $x_1\in(a/2,a)$ such that $h'_l(x_1)=h_r'(x_1)$. Since $\lim_{x\to a}h_l(x)<\lim_{x\to a}h_r(x)$ there must be at least $x_+\in(a/2,a)$ such that $h_l(x_+)=h_r(x_+)$. If there would be two different such points, for instance $x_+<x'_+$ then by the generalized mean value theorem there exists $\xi_+<\xi'_+$ such that
\begin{equation}\label{Mean_Value_Argument}
\begin{split}
1=\frac{h_r(x'_+)-h_r(x_+)}{h_l(x'_+)-h_l(x_+)}=\frac{h'_r(\xi'_+)}{h'_l(\xi'_+)}\hspace{0.5cm}\text{ and }\hspace{0.5cm}1=\frac{h_r(x_+)-h_r(a/2)}{h_l(x_+)-h_l(a/2)}=\frac{h'_r(\xi_+)}{h'_l(\xi_+)}
\end{split}
\end{equation}
in other words $h'_r(\xi'_+)=h'_l(\xi'_+)$ and $h'_r(\xi_+)=h'_l(\xi_+)$, a contradiction. Due to point symmetry $a/2$ then is a local maximum and $x_+$ as well as $x_-:=a-x_+$ are global minima of the free energy as a function of the first variable $\nu_1$ on $D^k_{\bar\nu}$. By Remark \ref{Permutation_Invariance}, $\nu_1\geq\nu_k$ and since $x_+>x_-$ we have $\nu_1=x_+$ and $\nu_k=x_-$. In particluar $\nu_1>\tilde u$ and thus by Lemma \ref{Only_Two_Poss} part (ii) the free energy minimizer has the form
\begin{equation*}
\bar\nu=(\underbrace{x_+,\dots,x_+}_{{}l\text{ times}},x_-,\dots,x_-)^T \text{ with }2\leq l<k.
\end{equation*}
Moreover if $l>1$, $\nu_1+\nu_l=2x_+>a>2\tilde u$ and hence by the same arguments as above $\nu_1>\nu_l$, a contradiction.
\begin{figure}[h]
\begin{center}
\includegraphics[width=4.5cm]{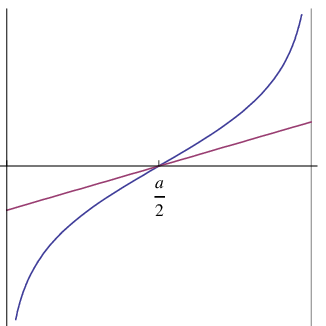}
\includegraphics[width=4.5cm]{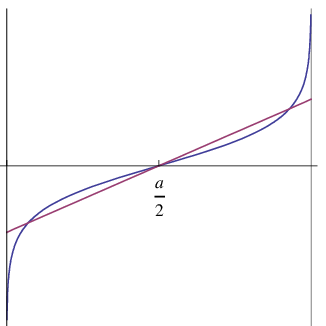}
\includegraphics[width=4.5cm]{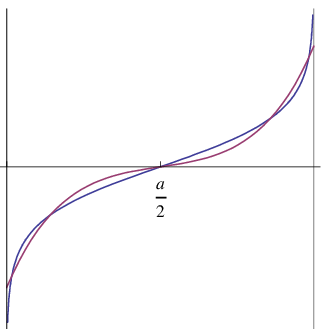}
\end{center}
\caption{\scriptsize{On the left side, $h_r$ and $h_l$ in the cases $2\leq z\leq4$ subcase one and $z>4$ subcase two. In the middle,  $h_r$ and $h_l$ in the cases $2\leq z\leq4$ subcase two and $z>4$ subcase one. On the right side, $h_r$ and $h_l$ in the case $z>4$ subcase four.}}
\label{Point_Symmetric_Functions}
\end{figure}

\bigskip
Now let $z>4$. We show, there is either one, two or no additional points $x\in(a/2,a]$ such that $h'_l(x)=h'_r(x)$. Let us look again at $\b_{z,a}$ defined in \eqref{Beta_For_h'}. For $z>4$, $\b_{z,a}$ has a local maximum in $a/2$ since $\b'_{z,a}(a/2)=0$ which can easily be seen from equation \eqref{Beta_Ableitung} and $\b''_{z,a}(a/2)=-(z-4)(\frac{a}{2})^{-(z+1)}<0$. We show, there is only one solution $\b_{z,q}'(x)=0$ on $(a/2,a]$ which must be a global minimizer since $\lim_{x\to a}\b_{z,a}(x)=\infty$. Indeed from \eqref{Beta_Ableitung_Replaced} we see, requiring $\b_{z,q}'$ to be zero is equivalent to the fixed point equation
\begin{equation*}
y=\frac{(z-y-1)y^{z-2}+1}{z-1}=:r_z(y)
\end{equation*}
having an unique solution on $[0,1)$. The r.h.s has the following properties: $r_z(0)=1/(z-1)>0$; $r_z(1)=1$; $r'_z(1)>1$ and $r_z$ is convex, since $r''_z(y)=(z^2-(3+y)z+2)y^{z-3}>0$. Combining these properties gives the uniqueness of the fixed point and thus the uniqueness of the extremal value of $\b_{z,a}$ which is a minimum that we want to call $\b_0(z,a)$.

\bigskip
Subcase one, let $a/2\geq\tilde u$. This is equivalent to $\b\geq\frac{1}{z-1}(\frac{a}{2})^{1-z}$ and hence by the exact same arguments as in the case $z\leq4$ subcase two, the free energy minimizer has the form $\bar\nu=(x_+,x_-,\dots,x_-)^T$ with $x_+>x_-$.

\bigskip
Subcase two, let $a/2<\tilde u$ and $\b<\b_0(z,a)$. Then we are in a situation as in 
case $z\leq4$ subcase one.  
In particular the free energy minimizer is the equidistribution.

\bigskip
Subcase three, let $a/2<\tilde u$ and $\b=\b_0(z,a)<\frac{1}{z-1}(\frac{a}{2})^{1-z}$. 
In this case, there is exactly one $x_1\in(a/2,a]$ such that $h_l'(x_1)=h_r'(x_1)$ and hence by the mean value argument already presented in \eqref{Mean_Value_Argument} there cannot be more then one $x_+\in(a/2,a]$ such that $h_l(x_+)=h_r(x_+)$. If no such $x_+$ exists, we are in the same situation as in subcase two (right above.) If such $x_+$ exists it must belong to a touching point of the graphes of $h_l$ and $h_r$ since otherwise because of $\lim_{x\to a}h_l(x)<\lim_{x\to a}h_r(x)$ there must be another point $(a/2,a]\ni\bar x_+\neq x_+$ with $h_l(\bar x_+)=h_r(\bar x_+)$. If it is a touching point of $h_l$ and $h_r$, then the free energy as a function of the first entry cannot attain a minimum in $x_+$, instead it is a saddle point and the minimum is attained in 
$a/2$. Consequently the minimizing distribution of the free energy is the equidistribution.

\bigskip
Subcase four, let $a/2<\tilde u$ and $\b_0(z,a)<\b<\frac{1}{z-1}(\frac{a}{2})^{1-z}$. In this case we have exactly two points $x_1<x_2$ such that $h'_l(x_i)=h'_r(x_i)$ with $i\in\{1,2\}$ and again by the mean value argument \eqref{Mean_Value_Argument} there cannot be more than two points $x_+>x'_+$ with $h_l(x_+)=h_r(x_+)$ and $h_l(x'_+)=h_r(x'_+)$. If no such point or only one such point exists, we can apply the same arguments as in subcase three (right above) and the equidistribution is the free energy minimizer. If both points exist and both belong to touching points of the graphes of $h_l$ and $h_r$ then again the equidistribution must be the minimizer. The case that both points exist and only one is a touching point is impossible. 

Now if both points exist and belong to real intersections of the graphes of $h_l$ and $h_r$, then we have three local minima attained in $x_-<a/2<x_+$ with $x_-:=a-x_+$.
Hence for $\nu_1$ the local minimizers $a/2$ and $x_+$ are competing to be the global minimizers. If $a/2$ is the global minimizer then by  Lemma \ref{Only_Two_Poss} the free energy is minimized by the equidistribution. If $x_+$ is the global minimizer, then notice if $x_+\leq\tilde u$ again by Lemma \ref{Only_Two_Poss} $x_+=a/2$ which contradicts $x_+>a/2$. Hence $x_+>\tilde u$ and the free energy minimizer has the form
\begin{equation*}
\bar\nu=(\underbrace{x_+,\dots,x_+}_{{}l\text{ times}},x_-,\dots,x_-)^T \text{ with }2\leq l<k.
\end{equation*}
Moreover if $l>1$, $\nu_1+\nu_l=2x_+>2\tilde u$ and hence by subcase one $\nu_1>\nu_l$, a contradiction.

\medskip
Finally, in order to have the minimizers in the format given in the theorem, define $u\in[0,1)$ such that $(1-u)/q=x_-$. This is always possible since $0<x_-\leq 1/q$. Of course $x_+=(1+(q-1)u)/q$.
$\Cox$

\bigskip
For the proof of part (ii) of \ref{Minimizer_Theorem} we need the following lemmata.
\begin{lem}\label{Key_Lemma_Minimizer}
For $q>2$ and $z\geq2$ there exist two temperatures $0<\b_0(q,z)<\b_1(q,z)$ such that for $0<\b<\b_0$ the mean-field equation only has the trivial solution $u=0$. For $\b_0<\b<\b_1$ the mean-field equation has two additional solutions $0<u_1<u_2<1$. Finally for $\b=\b_0$ or $\b\geq\b_1$ there is only one additional solution $0<u_2<1$.
\end{lem}

\textbf{Proof: }Let us write the temperature as a function of positive solutions of the mean-field equation 
\begin{equation}\label{Beta_Funktion}
\b_{q,z}(u):=q^{z-1}\frac{\log(1+(q-1)u)-\log(1-u)}{(1+(q-1)u)^{z-1}-(1-u)^{z-1}}
\end{equation}
\begin{figure}[h]
\begin{center}
\includegraphics[width=7cm]{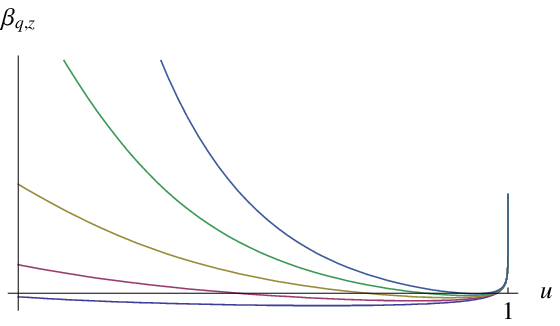}
\includegraphics[width=7cm]{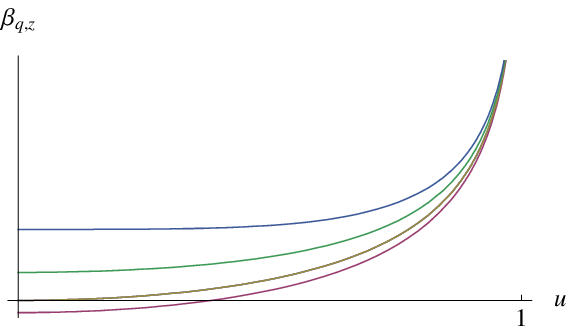}
\end{center}
\caption{\scriptsize{On the left side, $\b_{q,z}$ for $q=3$ and $z=3,\dots,7$. The cup shape of the graphs is a common feature for the parameter regimes $(q,z)\in([2,\infty)\times[2,\infty])\setminus(\{2\}\times[2,4])$. On the right side, $\b_{q,z}$ for $q=2$ and $z=2,\dots,4$. Here $\b_{q,z}$ is strictly increasing and this is a common feature for the parameter regimes $(q,z)\in\{2\}\times[2,4]$.}}
\label{Temperature_Function}
\end{figure}
and define $\lim_{u\to0}\b_{q,z}(u)=\frac{q^{z-1}}{z-1}=:\b_1$. Notice $\lim_{u\to0}\b_{q,z}'(u)=-\frac{1}{2}(q-2)q^{z-1}<0$ and $\lim_{u\to1}\b_{q,z}(u)=\infty=\lim_{u\to1}\b_{q,z}'(u)$. We will show, that $\b_{q,z}$ has exactly one extremal point attained in $u_0\in(0,1)$. This must be a local and hence global minimum that we want to call $\b_0$. Let us calculate
\begin{equation*}
\begin{split}
&0=\b_{q,z}'(u)=q^{z-1}\Bigl(\frac{q}{(1+(q-1)u)(1-u)[(1+(q-1)u)^{z-1}-(1-u)^{z-1}]}-\cr
&\frac{[\log(1+(q-1)u)-\log(1-u))](z-1)[(1+(q-1)u)^{z-2}(q-1)+(1-u)^{z-2}]}{[(1+(q-1)u)^{z-1}-(1-u)^{z-1}]^2}\Bigr).
\end{split}
\end{equation*}
Replacing $v:=(1+\frac{qu}{1-u})^{z-1}$ we can write equivalently
\begin{equation*}
\begin{split}
q-1&=v^{\frac{1}{z-1}}\frac{\log v-v+1}{v-1-v\log v}=:F_z(v).
\end{split}
\end{equation*}
Notice $u\mapsto v$, $(0,1)\mapsto(1,\infty)$ is strictly increasing and bijective. It suffices to show, that $F_z$ is bijective on $F_z^{-1}(1,\infty)$. First we have $\lim_{v\to1}F_z(v)=1$ and $\lim_{v\to\infty}F_z(v)=\infty$. We show that $F_z$ is strictly increasing on $F_z^{-1}(1,\infty)$ and calculate
\begin{equation*}
\begin{split}
0=F_z'(v)=v^{\frac{2-z}{z-1}}\frac{(z-2)v\log^2 v+(v^2-1)\log v -(v-1)^2}{(z-1)(1-v+v\log v)^2}
\end{split}
\end{equation*}
which is equivalent to 
\begin{equation*}
\begin{split}
z=\frac{(v^2-1)\log v-2v\log^2 v}{(v-1)^2-v\log^2 v}=:G(v).
\end{split}
\end{equation*}
Since $G(v)>4$ on $(1,\infty)$ (which we will see right below) for $2\leq z\leq4$ there are no extremal points of $F_z$ and in particular $F_z$ is bijective on $F_z^{-1}(1,\infty)$. Since $G(v)$ is also strictly increasing (which we will also see right below) and for $z>4$, $\lim_{v\to1}F'_z(v)=\frac{4-z}{3(z-1)}<0$ there is exactly one extremal points of $F_z$ which must be a minimum. In particular that minimum is smaller than one and hence $F_z$ is bijective on $F_z^{-1}(1,\infty)$.

To see that $G(v)>4$ and strictly increasing, use $\lim_{v\to1}G(v)=4$ and show $0<G'$ which is equivalent to 
\begin{equation*}
\begin{split}
G_1(v):=(v-1)^3(v+1)-6v(v-1)^2\log v+3v(v^2-1)\log^2v-v(v^2+1)\log^3 v>0.
\end{split}
\end{equation*}
One way to see that this is true is to show strict convexity of $G_1$ and use $\lim_{v\to1}G_1(v)=\lim_{v\to1}G_1'(v)=0$ and $\lim_{v\to\infty}G_1(v)=\infty$. Here $G_1''>0$ is equivalent to 
\begin{equation*}
\begin{split}
G_2(v):=4(v-1)^3-4(v-1)^2\log v+(v^2-1)\log^2 v-2v^2\log^3v>0
\end{split}
\end{equation*}
and again $\lim_{v\to1}G_2(v)=\lim_{v\to1}G_2'(v)=0$ and $\lim_{v\to\infty}G_2(v)=\infty$. Now again we want to show strict convexity of $G_2$, but this is equivalent to 
\begin{equation*}
\begin{split}
G_3(v):=1+4v-17v^2+12v^3+\log v-7v^2\log v-8v^2\log^2v-2v^2\log^3 v>0
\end{split}
\end{equation*}
and as before $\lim_{v\to1}G_3(v)=\lim_{v\to1}G_3'(v)=0$ and $\lim_{v\to\infty}G_3(v)=\infty$. Now again we want to show strict convexity of $G_3$, but this is equivalent to 
\begin{equation*}
\begin{split}
G_4(v):=-1-71v^2+72v^3-74v^2\log v-34v^2\log^2v-4v^2\log^3v>0
\end{split}
\end{equation*}
and $\lim_{v\to1}G_4(v)=\lim_{v\to1}G_4'(v)=0$ and $\lim_{v\to\infty}G_4(v)=\infty$. Now as above we want to show strict convexity of $G_4$, but this is equivalent to 
\begin{equation*}
\begin{split}
G_5(v):=54(v-1)-47\log v-13\log^2v-\log^3v>0
\end{split}
\end{equation*}
and now $\lim_{v\to1}G_5(v)=0$, $\lim_{v\to1}G_4'(v)=7$ and $\lim_{v\to\infty}G_4(v)=\infty$. Finally the strict convexity of $G_5$ is equivalent to $21+20\log v+3\log^2v>0$. But this is true and hence the above cascade gives $0<G'$. This finishes the proof.
$\Cox$

\begin{lem}\label{Key_Lemma_Minimizer_Curie_Weiss_High_Exp}
For $q=2$ and $z>4$ there exist two temperatures $0<\b_0(2,z)<\b_1(2,z)$ such that for $0<\b<\b_0$ the mean-field equation only has the trivial solution $u=0$. For $\b_0<\b<\b_1$ the mean-field equation has two additional solutions $0<u_1<u_2<1$. Finally for $\b=\b_0$ or $\b\geq\b_1$ there is only one additional solution $0<u_2<1$.
\end{lem}

\textbf{Proof: }$\b_{2,z}$ as defined in \eqref{Beta_Funktion} has the following properties: $\lim_{u\to0}\b_{2,z}'(u)=0$; $\lim_{u\to0}\b_{2,z}''(u)=2^{z-1}(4-z)/3<0$ and $\lim_{u\to1}\b_{2,z}(u)=\infty=\lim_{u\to1}\b_{2,z}'(u)$. Define $\lim_{u\to0}\b_{2,z}(u)=\frac{2^{z-1}}{z-1}=:\b_1$. Using the exact same arguments as presented in the proof of Lemma \ref{Key_Lemma_Minimizer} one can again show that $\b_{2,z}$ has exactly one extremal point $\b_0$ attained in $u_0\in(0,1)$. As before, the indicated parameter regimes are an immediate consequence of this fact.
$\Cox$

\begin{lem}\label{Key_Lemma_Minimizer_Curie_Weiss_Low_Exp}
For $q=2$ and $2\leq z\leq4$ there exist only one temperature $0<\b_1(2,z)$ such that for $0<\b\leq\b_1$ the mean-field equation only has the trivial solution $u=0$. For $\b>\b_1$ there is one additional solution $0<u_1<1$.
\end{lem}
\textbf{Proof: }$\b_{2,z}$ as defined in \eqref{Beta_Funktion} has the following properties: $\lim_{u\to0}\b_{2,z}'(u)=0$, $\lim_{u\to1}\b_{2,z}(u)=\infty=\lim_{u\to1}\b_{2,z}'(u)$; $\lim_{u\to0}\b_{2,z}''(u)=2^{z-1}(4-z)/3>0$ for $z<4$ and $\lim_{u\to0}\b_{2,4}''(u)=2^{z-1}(4-z)/3=0$; $\lim_{u\to0}\b_{2,4}'''(u)=0$; $\lim_{u\to0}\b_{2,4}''''(u)=64/5>0$. As a consequence for $2\leq z\leq4$, $\b_{2,z}$ has a local minimum in zero. We show, $\b_{2,z}$ is strictly increasing. Indeed $\b_{2,z}'>0$ is equivalent to 
\begin{equation*}
\begin{split}
F_z(v):=(v-1)(v^\frac{1}{z-1}+1)-(v+v^\frac{1}{z-1})\log v>0
\end{split}
\end{equation*}
for $v\in(1,\infty)$ where we made the one-to-one replacement $v=(\frac{1+u}{1-u})^{z-1}$. Notice $z\mapsto F_z$ is strictly decreasing point wise  since $d/dz F_z(v)<0$ is equivalent to $\log v<v-1$ which is of course true for all $v\in(1,\infty)$. Now in order to show $F_4>0$ we again use a cascade of convex functions. First, $F_4(1)=0$, $F'_4(1)=0$ and $F''_4>0$ is equivalent to $G(v):=5-9v^{2/3}+4v+2\log(v)>0$. Second, $G(1)=0$, $G'(1)=0$ and $G''>0$ is equivalent to $v>1$, but this is true. 

Consequently $\b_1(2,z):=\lim_{u\to0}\b_{2,z}(u)=2^{z-1}/(z-1)$.
$\Cox$

\bigskip
\textbf{Proof of Theorem \ref{Minimizer_Theorem} part (ii): }
The above lemmata consider the temperature parameter as a function of positive solutions of the mean-field equation 
\begin{equation*}
\b_{q,z}(u)=q^{z-1}\frac{\log(1+(q-1)u)-\log(1-u)}{(1+(q-1)u)^{z-1}-(1-u)^{z-1}}.
\end{equation*}
This function is positive. 

\medskip
In the parameter regimes considered in Lemma \ref{Key_Lemma_Minimizer} and Lemma \ref{Key_Lemma_Minimizer_Curie_Weiss_High_Exp} $\b_0$ is the unique global minimum of $\b_{q,z}$ and $\b_1=\lim_{u\to0}\b_{q,z}(u)=\frac{q^{z-1}}{z-1}$. Let us connect this with the free energy as a function of $u$.
\begin{equation}\label{One_Input_Free_Energy}
\begin{split}
\G_{\b,q,z}(\bar{\nu})=&-\frac{\b}{z}\bar{\nu}_1^{z}+\bar{\nu}_1\log(q\bar{\nu}_1)+
(q-1)(-\frac{\b}{z}\bar{\nu}_2^{z}+\bar{\nu}_2\log(q\bar{\nu}_2))\cr
=&\frac{1}{q}[(1+(q-1)u)\log(1+(q-1)u)+(q-1)(1-u)\log(1-u)]\cr
&-\frac{\b}{z}q^{-z}[(1+(q-1)u)^{z}+(q-1)(1-u)^{z}]=:k_{\b,q,z}(u)
\end{split}
\end{equation}
\begin{figure}[h]
\begin{center}
\includegraphics[width=7cm]{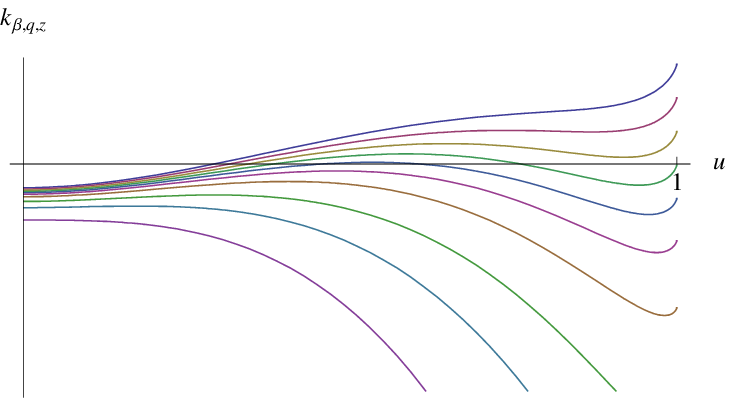}
\includegraphics[width=7cm]{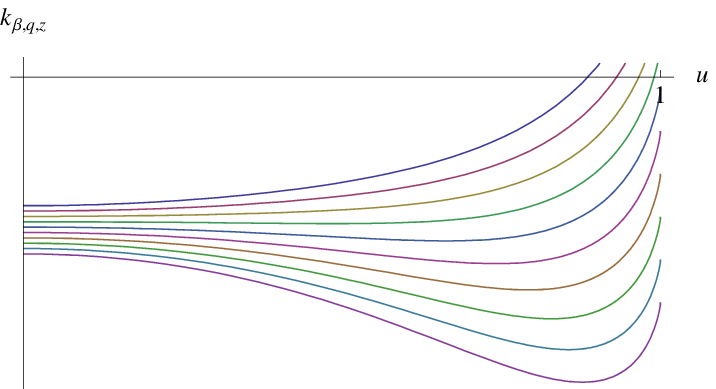}
\end{center}
\caption{\scriptsize{On the left side, $k_{\b,q,z}$ for $q=3$, $z=4$ and $\b=3.8,\dots,9$. The fact that the shape of the graph changes from a single local minimum attained away from zero, to two local minima and back again to one locals minimum attained in zero is a common feature in the parameter regimes $(q,z)\in([2,\infty)\times[2,\infty])\setminus(\{2\}\times[2,4])$. One can clearly see the first-order nature of the phase-transition. On the right side, $k_{\b,q,z}$ for $q=2$, $z=4$ and $\b=2.4,2.5,\dots,3.3$. The fact that the point where the global minimum is attained moves into zero from the right is a common feature for the parameter regimes $(q,z)\in\{2\}\times[2,4]$. This indicates a second-order phase-transition.}}
\label{Free_Energy}
\end{figure}
and its derivatives
\begin{equation}
\begin{split}
k_{\b,q,z}'(u)=&-\frac{q-1}{q^{z}}\b[(1+(q-1)u)^{z-1}-(1-u)^{z-1}]-\frac{q-1}{q}\log\frac{1-u}{1+(q-1)u},\cr
k_{\b,q,z}''(u)=&-\frac{q-1}{q^{z}}\b(z-1)[(q-1)(1+(q-1)u)^{z-2)}+(1-u)^{z-2)}]\cr
&+\frac{q-1}{(1-u)(1+(q-1)u)}.
\end{split}
\end{equation}
Notice $k_{\b,q,z}'(0)=0$ and $k$ has a local minimum in zero iff $\b<\frac{q^{z-1}}{z-1}=\b_1$.
Since also $\lim_{u\to1}k_{\b,q,z}'(u)=+\infty$ we can assert the following:

\begin{enumerate}
\item If $\b<\b_0<\b_1$ then in $u=0$ the free energy must attain its global minimum.
\item If $\b\geq\b_1$ then in zero there is a local maximum and according to Lemma \ref{Key_Lemma_Minimizer} and Lemma \ref{Key_Lemma_Minimizer_Curie_Weiss_High_Exp} there is exactly one more extremal point, but this must be a global minimum.
\item If $\b=\b_0<\b_1$ the additional extremal point must be a saddle point since if it would be a local maximum, then there must be another local minimum and hence another extremal point, but the additional extremal point is the only one.
\item If $\b_0<\b<\b_1$ then the two additional extremal points $u_1<u_2$ are either two saddle points or a local maximum (attained in $u_1$) and a local minimum (attained in $u_2$.)
\end{enumerate}

Since $\frac{d}{d\b}k_{\b,q,z}(u)=-\frac{q^{-z}}{z}[(1+(q-1)u)^{z}+(q-1)(1-u)^{z}]<0$ the free energy decreases for every $u$ if $\b$ increases. Since $\frac{d}{du}\frac{d}{d\b}k_{\b,q,z}(u)=-\frac{q-1}{q^{z}}[(1+(q-1)u)^{z-1}-(1-u)^{z-1}]<0$, for larger $u$ this decrease is also strictly larger and hence for $\b$ moving up from $\b_0$ to $\b_1$, $k_{\b,q,z}(u_2)$ is going down faster than $k_{\b,q,z}(0)$. Since for $\b\geq\b_1$, $u_2$ becomes the global minimum, and $k$ is continuous w.r.t every parameter, there must be a $\b_0<\b_c\leq\b_1$ where $k_{\b,q,z}(0)=k_{\b,q,z}(u_2)$ and indeed for $\b>\b_c$ the minimizer of the free energy $\G_{\b,q,z}$ is defined by the largest solution of the mean-field equation. 

\medskip
In the parameter regime considered in Lemma \ref{Key_Lemma_Minimizer_Curie_Weiss_Low_Exp} the situation is simpler and we can set $\b_0=\b_1=\b_c$. In particular
\begin{enumerate}
\item If $\b<\b_c$ then in $u=0$ the free energy must attain its global minimum.
\item If $\b>\b_c$ then in zero there is a local maximum and according to Lemma \ref{Key_Lemma_Minimizer_Curie_Weiss_Low_Exp} there is exactly one more extremal point, but this must be a global minimum.
\end{enumerate}
$\Cox$

\bigskip
\textbf{Proof of Theorem \ref{Minimizer_Theorem} part (iii): }In the cases $z\geq2$, $q\geq2$ and $z>4$, $q=2$ we have $\b_0<\b_c\leq\b_1$ and $\lim_{\b\searrow\b_c}u(\b,q,z)=u_2(q,z)>0=\lim_{\b\nearrow\b_c}u(\b,q,z)$ where we used notation from the proof of part 2 of \ref{Minimizer_Theorem} with $u_2(q,z)=u_2$. Hence $\b\mapsto u(\b,q,z)$ is discontinuous in $\b_c$.

In the case $2\leq z\leq4$, $q=2$ we have $\lim_{\b\searrow\b_c}u(\b,q,z)=0=\lim_{\b\nearrow\b_c}u(\b,q,z)$ by the monotonicity of $u\mapsto\b_{q,z}(u)$ and hence $\b\mapsto u(\b,q,z)$ is continuous in $\b_c$.
$\Cox$


\subsection{Proof of Proposition \ref{q_Monotonicity_Of_Beta}}
It suffice to show $\partial_q\b_c(q,z)\geq0$, where $\partial_q\b_c$ stands for the partial derivative of $\b_c$ in the direction $q$. Without restriction we consider $3\leq q\in\R$. We know that $\b_c>0$ and the corresponding value $u_c\in(0,1)$ are solutions of the equations:
\begin{equation}\label{ku=knull}
k_{\b,q,z}(u)=k_{\b,q,z}(0)=-\frac{\b}{z}q^{1-z}\hspace{0.5cm}\text{ and }\hspace{0.5cm}k_{\b,q,z}'(u)=0,
\end{equation}
where $k_{\b,q,z}$ is given in \eqref{One_Input_Free_Energy}. The first condition is equivalent to
\begin{equation}\label{F}
\begin{split}
&F(\b,q,u):=\b f(q,u)+g(q,u):=\cr
&-\frac{\b}{z}q^{-z}[(1+(q-1)u)^z+(q-1)(1-u)^z-q]\cr
&+\frac{1}{q}[(1+(q-1)u)\log(1+(q-1)u)+(q-1)(1-u)\log(1-u)]=0.
\end{split}
\end{equation}
The second condition is equivalent to 
\begin{equation*}\label{G}
\begin{split}
G(\b,q,u):=\partial_u F(\b,q,u)=:\b\partial_u f(q,u)+\partial_u g(q,u)=0.
\end{split}
\end{equation*}
Taking the derivative along a path of solutions we get a two-dimensional system of equations
\begin{align*}
\frac{d}{dq}F(\b(q),q,u(q))=\partial_\b F(\b,q,u)\partial_q\b(q)+\partial_q F(\b,q,u)+\partial_u F(\b,q,u)\partial_q u(q)&=0\cr
\frac{d}{dq}G(\b(q),q,u(q))=\partial_\b G(\b,q,u)\partial_q\b(q)+\partial_q G(\b,q,u)+\partial_u G(\b,q,u)\partial_q u(q)&=0,
\end{align*}
where we wrote for simplicity $\b_c(q)=\b_c(q,z)$. This is equivalent to
\begin{align*}
\begin{pmatrix}
\partial_q\beta(q)\\
\partial_q u(q)
\end{pmatrix}
=-\begin{pmatrix}
\partial_\beta F(\b,q,u)&\partial_u F(\b,q,u)\\
\partial_\beta G(\b,q,u)&\partial_u G(\b,q,u)
\end{pmatrix}^{-1}
\begin{pmatrix}
\partial_q F(\b,q,u)\\
\partial_q G(\b,q,u)
\end{pmatrix}
\end{align*}
which leads to
\begin{align*}\label{dbq}
\partial_q\beta(q)=-\frac{\partial_u G(\b,q,u)\partial_q F(\b,q,u)-\partial_u F(\b,q,u)\partial_q G(\b,q,u)}{\partial_\beta F(\b,q,u)\partial_u G(\b,q,u)-\partial_\beta G(\b,q,u)\partial_u F(\b,q,u)}.
\end{align*}
Now we can use that for our solutions $G(\b,q,u)=\partial_u F(\b,q,u)=0$ and thus we have
\begin{equation*}\label{dbq_simp}
\partial_q\b_c(q)=-\frac{\partial_q F(\b,q,u)}{\partial_\b F(\b,q,u)}=-\frac{\partial_q F(\b,q,u)}{f(q,u)}.
\end{equation*}
Notice $f(q,u)<0$ since $f(q,0)=0$ and 
\begin{equation*}
\partial_u f(q,u)=q^{-z}(q-1)[(1-u)^{z-1}-(1+(q-1)u)^{z-1}]<0.
\end{equation*}
where we used $1-u<1+(q-1)u$. Hence it suffices to show 
\begin{equation}\label{Fq}
\partial_q F(\b,q,u)=\b \partial_q f(q,u)+\partial_q g(q,u)\geq 0. 
\end{equation}
A solution of \eqref{F} satisfies $\b=-g(q,u)/f(q,u)$. Thus we can eliminate $\b$ in (\ref{Fq}) and show instead
\begin{equation}\label{inequ}
\partial_q f(q,u)g(q,u)-\partial_qg(q,u)f(q,u)\geq 0.
\end{equation}
It would be sufficient to show that \eqref{inequ} is true for solutions of \eqref{ku=knull}. Nevertheless, we will prove \eqref{inequ} for all $q\in\R_{+}$, $z\in\R_{+}$ and $u\in[0,1]$. Multiplying \eqref{inequ} with $zq^{z+2}$, the inequality becomes
\begin{align}\label{inequ_bar}
0&\leq zq^{z+1}\partial_q f(q,u)\cdot qg(q,u)-q^2\partial_qg(q,u)\cdot zq^{z}f(q,u)\cr
&=\bar{f}_q(q,u)\cdot\tilde{g}(q,u)+\bar{g}_q(q,u)\cdot\tilde{f}(q,u),
\end{align}
with
\begin{align*}
\bar{f}_q(q,u):=&zq^{z+1}\partial_q f(q,u)
=z(1-u)(1+(q-1)u)^{z-1}+(q(z-1)-z)(1-u)^{z}-q(z-1)\cr
\tilde{g}(q,u):=&qg(q,u)=(1+(q-1)u)\log(1+(q-1)u)+(q-1)(1-u)\log(1-u)\cr
\bar{g}_q(q,u):=&q^2\partial_qg(q,u)=qu-(1-u)[\log(1+(q-1)u)-\log(1-u)]\cr
\tilde{f}(q,u):=&-zq^{z}f(q,u)=(1+(q-1)u)^{z}+(q-1)(1-u)^{z}-q.
\end{align*}
We have the following properties:

\begin{enumerate}
\item $u\mapsto\tilde{f}(q,u)\geq 0$ since $\tilde{f}(q,0)=0$ and $\partial_u\tilde{f}(q,u)=z(q-1)[(1+(q-1)u)^{z-1}-(1-u)^{z-1}]\geq 0$
\item $u\mapsto\tilde{g}(q,u)\geq 0$ since $\tilde{g}(q,0)=0$ and $\partial_u\tilde{g}(q,u)=(q-1)[\log(1+(q-1)u)-\log(1-u)]\geq 0$
\item $u\mapsto\bar{g}_q(q,u)\geq 0$ since $\bar{g}_q(q,0)=0$ and $\partial_u\bar{g}_q(q,u)=q-\frac{(q-1)(1-u)}{1+(q-1)u}+\log(1+(q-1)u)-\log(1-u)-1\geq 0$ since $q-1-(q-1)(1-u)/(1+(q-1)u)=q-q/(1+(q-1)u)>0$
\end{enumerate}
The more involved function is $u\mapsto\bar{f}_q(q,u)$ since it can be positive and negative. For the problematic case we define a set of $u$'s where $\bar{f}_q(q,u)$ is negative, i.e. $[0, 1]\supset A_q:=\{u\in[0, 1]:\bar{f}_q(q,u)<0\}$. Of course \eqref{inequ_bar} is true on $[0, 1]\setminus A_q$. Hence we only have to show on $A_q$ the inequality
\begin{equation*}
0\leq\bar{f}_q(q,u)\frac{\tilde{g}(q,u)}{\bar{g}_q(q,u)}+\tilde{f}(q,u).
\end{equation*}
Notice, $\bar{g}_q(u, q)=0$ only for $u=0$, but $0\notin A_q$ since $\bar{f}_q(0, q)=0$. We eliminate the fraction by the estimate
\begin{equation*}
\frac{\tilde{g}(u, q)}{\bar{g}_q(u, q)}\leq (q-1).
\end{equation*}
To see that this is true we use the following equivalent expressions:
\begin{equation*}
\begin{split}
\tilde{g}(q,u)&\leq (q-1)\bar{g}_q(q,u)\cr
(1+(q-1)u)\log(1+(q-1)u)&\leq(q-1)[qu-(1-u)\log(1+(q-1)u)]\cr
\log(1+(q-1)u)&\leq(q-1)u
\end{split}
\end{equation*}
Since $\bar{f}_q(q,u)$ is negative on $A_q$, we have $\bar{f}_q(q,u)(q-1)\leq \bar{f}_q(q,u)\frac{\tilde{g}(q,u)}{\bar{g}_q(q,u)}$ and all that is left to prove is
\begin{equation*}
0\leq\bar{f}_q(q,u)(q-1)+\tilde{f}(q,u).
\end{equation*}
Since $\bar{f}_q(q,0)(q-1)+\tilde{f}(q,0)=0$, it suffices to show
\begin{equation}\label{Final_Inequality}
\frac{d}{du}(\bar{f}_q(q,u)(q-1)+\tilde{f}(q,u))\geq 0.
\end{equation}
For simplicity let us write $A:=1-u$ and $B:=1+(q-1)u$, then \eqref{Final_Inequality} is true since the last of following equivalent expressions is clearly true
\begin{equation*}
\begin{split}
\partial_u (\bar{f}_q(q,u)(q-1)+\tilde{f}(q,u))&\geq 0\cr
(z-q(z-1))A^{z-1}-A^{z-1}+(q-1)(z-1)AB^{z-2}&\geq 0\cr
AB^{z-1}-A^{z-1}B&\geq 0\cr
\end{split}
\end{equation*}
$\Cox$

\section{Proof of Theorem \ref{Generalized_HaKu}}
Please note, most of the calculations done in this section work also for more general differentiable interaction functions. We prepare the proof by two propositions.
\begin{prop}
For each finite $N$ we have the representation 
\begin{equation}\label{Representation_Kernel}
Q^N_{\b,q,z,(r_1,\dots,r_s)}(k|\nu)= \frac{r_kA(\b_k,r_k,N_k)}{\sum_{l=1}^sr_lA(\b_l,r_l,N_l)}
\end{equation}
with $N_l=(N-1)\nu_l$, $\b_l=\b (N_l/N)^{z-1}$ and $A(\b,r,M)=\pi^{M}_{\b,r,z}(\exp(\b L_{M}^{\cdot}(1)^{z-1}))$.
\end{prop}
\textbf{Proof: }To compute the l.h.s of \eqref{Representation_Kernel} starting from the generalized fuzzy Potts measure, because of permutation invariance we can set $i=1$ and write for a fuzzy configuration $\eta$ on $\{2,\dots,N\}$
\begin{equation*}
\mu^N_{\b,q,z,(r_1,...,r_s)}(Y_1=k|Y_{\{2,\dots,N\}}=\eta)=\frac{1}{Z(\eta)}\sum_{\xi:T(\xi)=(k,\eta)}\pi^N_{\b,q,z}(\xi)
\end{equation*}
where $Z(\eta)$ is a normalization constant. Parallel to the proof of Proposition 5.2 in \cite{HK04}, it suffices to consider 
\begin{equation*}
\begin{split}
&\sum_{\xi:T(\xi)=(k,\eta)}\exp(\frac{\b N}{z}\sum_{i=1}^q(L_N^\xi(i))^z)=\sum_{\xi:T(\xi)=(k,\eta)}\exp(\frac{\b N}{zN^z}\sum_{i=1}^q(\sum_{j=1}^N1_{\xi_j=i})^z)\cr
&=\sum_{\xi:T(\xi)=(k,\eta)}[\exp(\frac{\b N}{zN^z}\sum_{i:T(i)=k}(1_{\xi_1=i}+\sum_{j\in\L_k}1_{\xi_j=i})^z)\prod_{l\neq k}\exp(\frac{\b N}{zN^z}\sum_{i:T(i)=l}(\sum_{j\in\L_l}1_{\xi_j=i})^z)]\cr
\end{split}
\end{equation*}
where we used $\L_l:=\{j\in\{2,\dots,N\}:\eta_j=l\}$.
Deviding this expression by $\prod_{l=1}^s\sum_{\xi_{\L_l}:T(\xi_{\L_l})=l}\exp(\frac{\b N}{zN^z}\sum_{i:T(i)=l}(\sum_{j\in\L_l}1_{\xi_j=i})^z)$ which is only dependent on $\eta$ gives
\begin{equation*}
\begin{split}
&\frac{\sum_{\xi_1:T(\xi_1)=k}\sum_{\xi_{\L_k}:T(\xi_{\L_k})=k}\exp(\frac{\b N}{zN^z}\sum_{i:T(i)=k}(1_{\xi_1=i}+\sum_{j\in\L_k}1_{\xi_j=i})^z)}{\sum_{\xi_{\L_k}:T(\xi_{\L_k})=k}\exp(\frac{\b N}{zN^z}\sum_{i:T(i)=k}(\sum_{j\in\L_k}1_{\xi_j=i})^z)}\cr
&=\sum_{\xi_1:T(\xi_1)=k}
\pi^{|\L_k|}_{\b,r_k,z}(\exp(\b N\sum_{i:T(i)=k}(\frac{|\L_k|}{N}L_{|\L_k|}^{\cdot}(i))^{z-1}\frac{1}{N}1_{\xi_1=i}+o(\frac{1}{N})))\cr
&=r_k\pi^{|\L_k|}_{\b,r_k,z}(\exp(\b(\frac{|\L_k|}{N}L_{|\L_k|}^{\cdot}(1))^{z-1}+o(1)))\cr
\end{split}
\end{equation*}
where we used Taylor expansion in the second last line. Since we are only interested in the limiting behavior of $Q^N$ as the system grows, by slight abuse of notation we can absorbe the asymptotic constant $o(1)$ into the normalization constant and hence the representation result follows.
$\Cox$

\begin{prop}\label{Limiting_Cond_Prob}
We have for boundary conditions $\nu^{(N)}\to\nu$, 
\begin{equation}
\label{limpropfp}
\lim_{N\to\infty}Q^N_{\b,q,z,(r_1,...,r_s)}(k|\nu^{(N)})=\frac{C(\b\nu_k^{z-1},r_k)}{\sum_{l=1}^s{C(\b\nu_l^{z-1},r_l)}}
\end{equation}
whenever $\nu_k^{z-1}\neq\b_c(r_k,z)/\b$ for all $r_k\geq2$ and $z\geq2$, where
\begin{align}
C(\b\nu_k^{z-1},r_k):=\notag
\begin{cases}
r_k\exp(\b(\frac{\nu_k}{r_k})^{z-1})&,\b\nu_k^{z-1}<\b_c(r_k,z)\cr
(r_k-1)\exp(\b\nu_k^{z-1}(\frac{1-u(\b\nu_k^{z-1},r_k,z)}{r_k})^{z-1})+&,\b\nu_k^{z-1}>\b_c(r_k,z)\cr
\hspace{0.2cm}\exp(\b\nu_k^{z-1}(\frac{(r_k-1)u(\b\nu_k^{z-1},r_k,z)+1}{r_k})^{z-1}).
\end{cases}\label{C_genfuzzy}
\end{align}
As a reminder, $u(\b\nu_k^{z-1},r_k,z)$ is the largest solution of the generalized mean-field equation \eqref{MFeqPhi}.
\end{prop}
\textbf{Proof: }The result is a direct consequence of the generalized Ellis-Wang Theorem \ref{Generalized_Ellis_Wang}.
$\Cox$

\bigskip
\textbf{Proof of Theorem \ref{Generalized_HaKu}: }
By Proposition \ref{Limiting_Cond_Prob}, for $2<z\leq4$ the points of discontinuity are precisely given by the values $\nu_k^{z-1}=\b_c(r_k,z)/\b$ for those $k\in\{1,\dots,s\}$ with $r_k\geq3$ for which $\b_c(r_k,z)/\b<1$. In particular if $r_i\leq2$ for all $i\in\{1,\dots,s\}$ no such points exist, this gives part (i). By Proposition \ref{q_Monotonicity_Of_Beta} $\b_c(r,z)$ is an increasing function of $r$, thus points of discontinuity can only be present if $\b$ is at least larger or equal than the critical inverse temperature of the smallest class that can have a second-order phase-transition. By picking two different approximating sequences of boundary conditions $\nu^{(N)}_k\searrow\nu_k$ and $\tilde\nu^{(N)}_k\nearrow\nu_k$ it is also clear that for those points of discontinuity the limit does not exist. This gives (ii) and (iii).
$\Cox$

\section{Appendix}
\subsection{Bifurcation analysis}
We have seen, that in different parameter regimes of the generalized Potts model different kinds of phase-transitions can appear. This is of course related to the apprearance (and disappearance) of local minima and maxima in the free energy as a function of $u\in[0,1]$ that we called $k_{\b,q,z}$ (see \eqref{One_Input_Free_Energy}.) A complete picture of possible bifurcations for general potentials is presented in \cite{PoSt78}. In this appendix we want to at least provide some figures showing the bifurcation phenomena that can appear in the generalized Potts model in particular. 
\begin{figure}[h]
\begin{center}
\includegraphics[width=3.5cm]{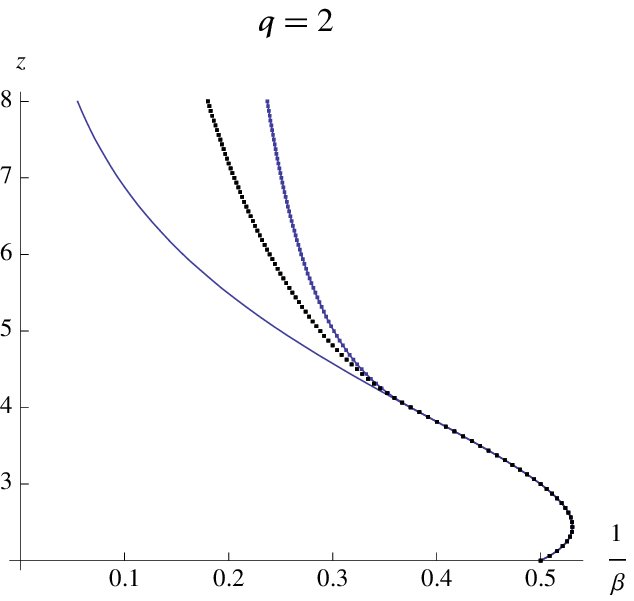}
\includegraphics[width=3.5cm]{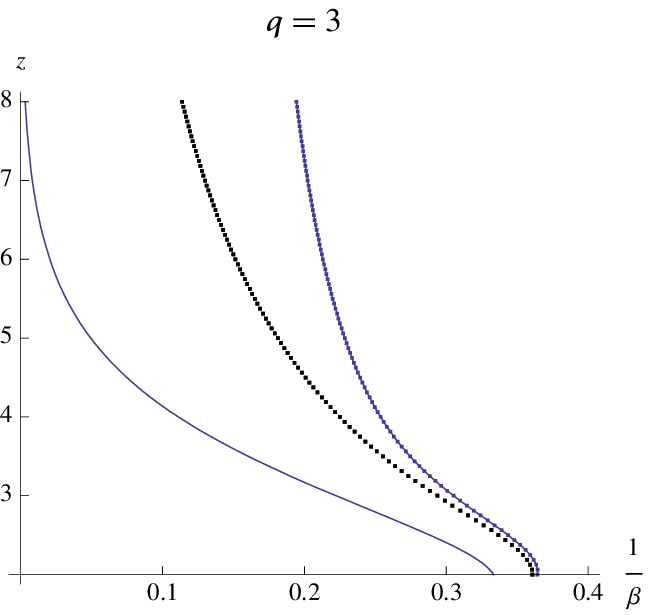}
\includegraphics[width=3.5cm]{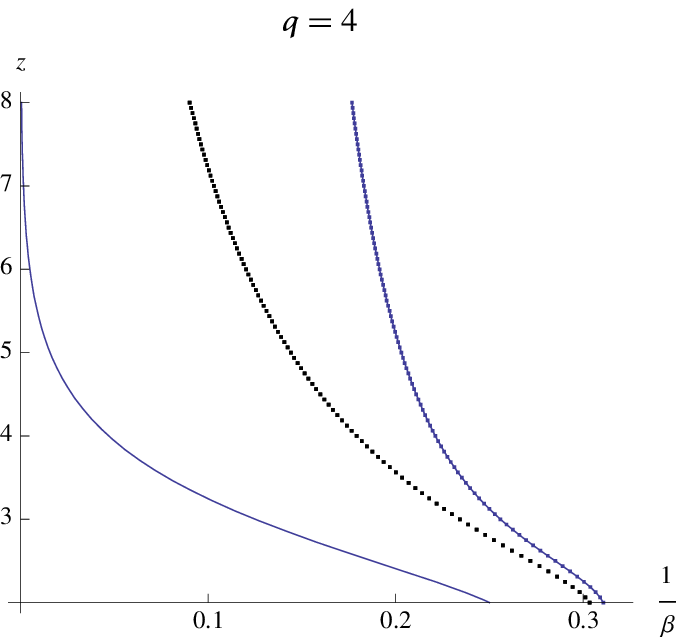}
\includegraphics[width=3.5cm]{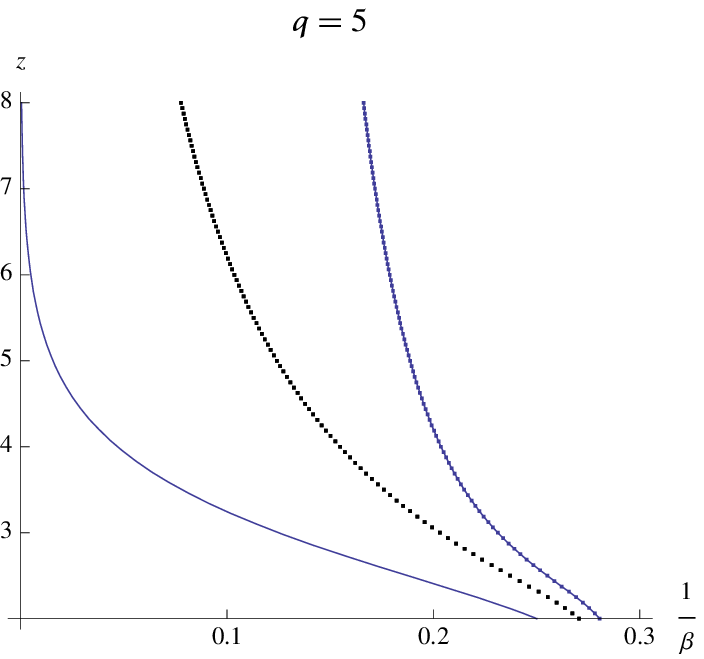}
\end{center}
\caption{\scriptsize{For $q=2$ the area left of the middle line is the phase-transition region. There is a triple point at $z=4$ where all extremal points fall in the same place, namely zero. Below $z=4$ there is the second-order phase-transition boundary and the three lines lie exactly on top of each other. Above $z=4$ there is a first-order phase-transition and the two additional lines right and left of the phase-transition boundary indicate bifurcation phenomena. To be more precise, the left line indicates where the local minimum at $u=0$ and the local maximum at $u_1\geq0$ join. The right line indicates where the local maximum at $u_1>0$ and the local minimum at $u_2\geq u_1$ join. Of course the phase-transition boundary must lie between these lines. We give a schematic picture for this in Figure \ref{Bifurcation_Schematisch}. For $q=3,4,5$ the situation is simpler since no second-order phase-transition is present.}}
\label{Bifurcation}
\end{figure}

\begin{figure}[h]
\begin{center}
\includegraphics[width=7cm]{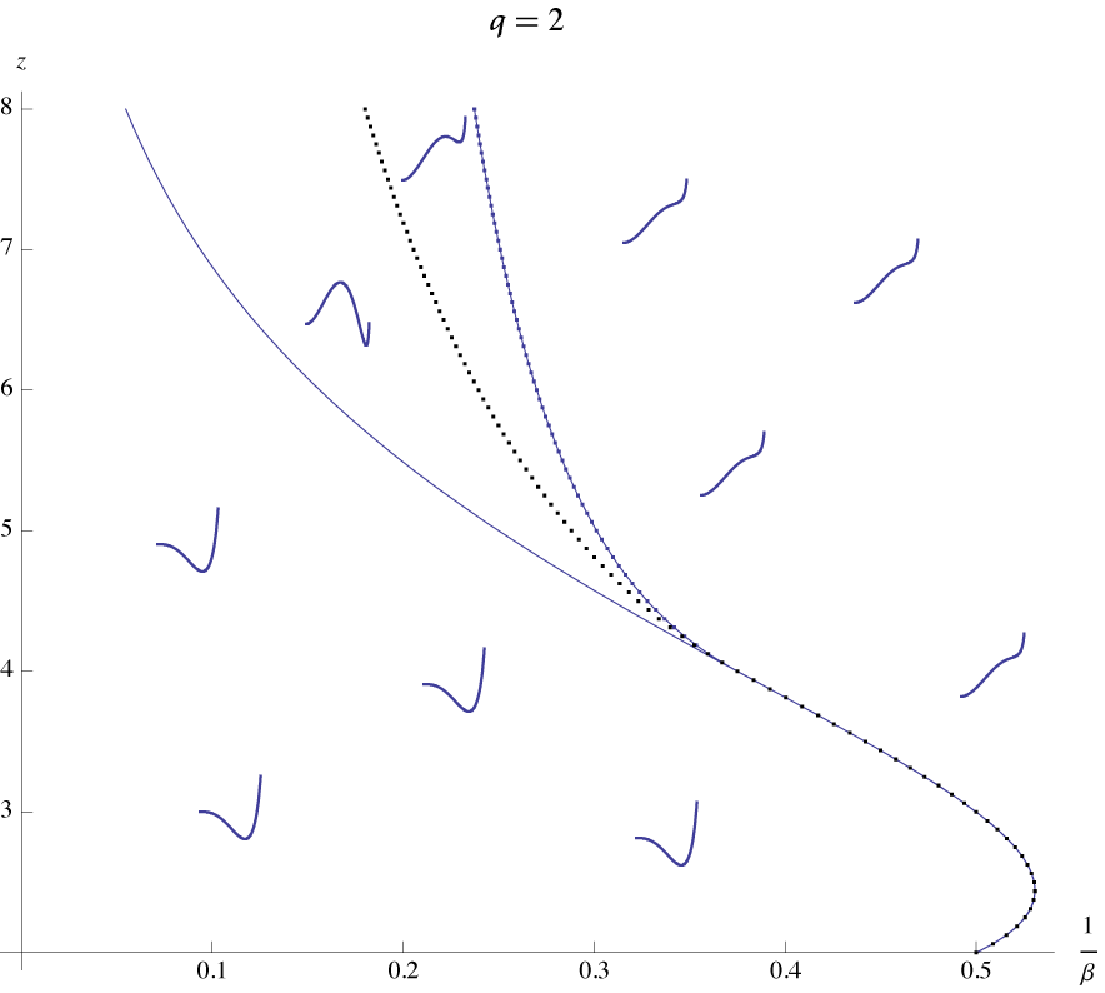}
\includegraphics[width=7cm]{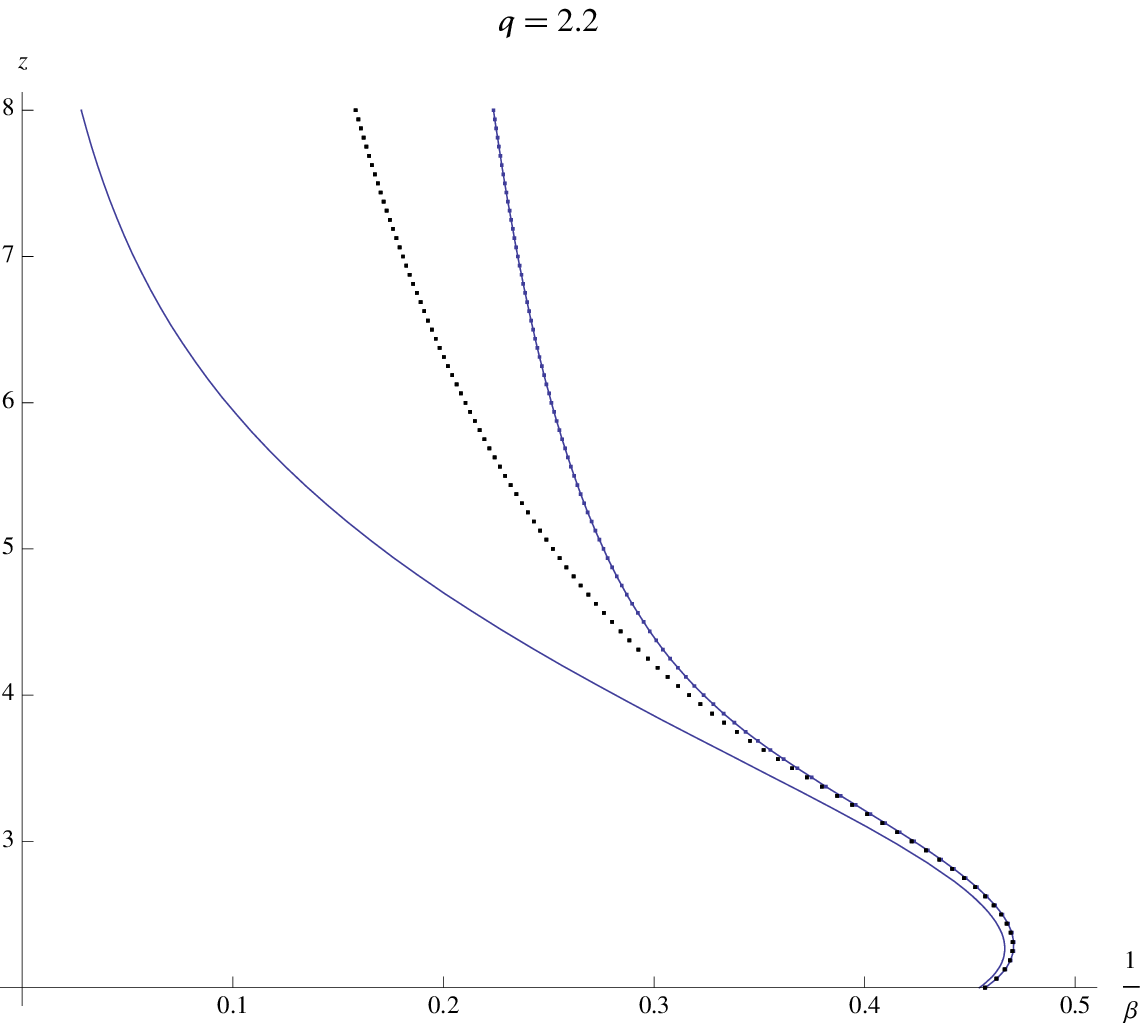}
\end{center}
\caption{\scriptsize{On the left: A schematic indication for the bifurcation phenomena present in case of $q=2$
where the small graphs are prototypical representations of the shape of the free energy. 
 On the right: For $q>2$ the bifurcation lines do not join and the phase-transition boundary lies in an area where there are two local minima and one local maximum present.}}
\label{Bifurcation_Schematisch}
\end{figure}

Note that only the left bifurcation line in each image in Figure \ref{Bifurcation} and Figure \ref{Bifurcation_Schematisch} we can compute exactly via $(\frac{d}{du})^2k_{\b,q,z}(u)|_{u=0}=0$ which is equivalent to $\frac{1}{\b}=\frac{z-1}{q^{z-1}}$. The right line in each of the same images shows $\b_0(q,z)$ as defined for example in Lemma \ref{Key_Lemma_Minimizer} which we computed numerically. The middle line showing $\b_c(q,z)$ we also calculated numerically.

\bigskip
On a computational level there is no reason not to assume $q$ to be continuous. In fact all our proofs work well with $\R\ni q\geq2$. We already showed that the possibility of a second-order phase-transition disappears for $q>2$. This one can also see in the bifurcation picture as indicated in Figure \ref{Bifurcation_Schematisch}.

\subsection{Random cluster representation and $z$-clique variables}\label{Random cluster representation and $z$-clique variables}

There is an equivalent notation for the Hamiltonian of the standard Potts model on the complete 
graph, namely $F_{\b,q,2}(L_N^{\xi})=-\frac{\b}{N^2}\sum_{1\leq i<j\leq N}1_{\xi_i=\xi_j}-\frac{\b}{N}$. For general integer valued exponents $z\geq2$ an equivalent notation for the Hamiltonian is given by
\begin{equation}\label{Alt_Hamiltonian}
F_{\b,q,z}(L_N^{\xi})=-\frac{\b (z-1)!}{N^z}\sum_ {D\subset\{1,\dots,N\},|D|=z}1_{\xi|_D=c}+O(\frac{1}{N})
\end{equation}
where $\xi|_D=c$ means that the given configuration $\xi$ has a constant $q$-coloring on the subset $D$ of size $z$. $N$ times the additional term is bounded by a constant as the system size grows and hence plays no role in the large deviation analysis and for the limiting Potts measure away from $\b_c(q,z)$. 

\medskip
We would like to describe now an extension of  the well-known random cluster 
representation of the nearest-neighbor Potts model on a general graph with $N$ vertices to  
interactions between $z=2,3,4,\dots $ spins. 
Denote by $\D$ a subset of the set of subsets of vertices $\{1,\dots,N\}$ with $z$ sites. In other words $\D$ is a subset of the $z$-cliques. This defines a graph in the usual sense when we say that there is an edge between sites $i\neq j$ iff there exists $D\in\D$ with $i,j\in D$.
We define the corresponding $\D$-Potts-Hamiltonian by 
\begin{equation}\label{Alt_Hamiltonian-zwei}
F_{\D}(\xi)=-\beta\sum_ {D\subset \D}1_{\xi|_D=c} 
\end{equation}
for a spin-configuration $\x\in \{1,\dots, N\}$. 
In the limit away from $\b_c(q,z)$ this corresponds to 
the generalized mean-field Potts measure for integer exponent $z$ when we take $\D$ to be 
the set of all subsets of $\{1,\dots,N\}$ with exactly $z$ elements. 

Let us now describe a random cluster representation for Gibbs measure 
corresponding to \eqref{Alt_Hamiltonian-zwei}.
Given $\D$, define the probability measure on $\{1,\dots, q\}^N \times \{0,1\}^\D$ by 
\begin{equation}\label{Coupling}
K(\sigma,\o)=C \prod_{D \in \D}
\Bigl((1-p)1_{\o(D)=0}+ p 1_{\o(D)=1}1_D(\s)\Bigr) 
\end{equation}
with $1_D(\s)$ the indicator of the event that $\s$ is constant on $D$ and $C$ the normalization. For $z=2$ this is the so-called \textit{Edwards-Sokal measure} presented in \cite{ES88}. Summing over the "clique-variables" $\o$ we get the marginal distribution on $\{1,\dots, q\}^N$  
\begin{equation*}
\begin{split}
\sum_{\o} K(\sigma,\o)&=C \sum_{\o}\prod_{D \in \D}\Bigl((1-p)1_{\o(D)=0}+ p1_{\o(D)=1}1_D(\s)\Bigr)\cr
&=C\prod_{D \in \D}\Bigl((1-p)+ p 1_D(\s)\Bigr)=C\prod_{D \in \D}(1-p)^{1-1_D(\s)}.
\end{split}
\end{equation*}
This equals the generalized Potts measure with Hamiltonian 
\eqref{Alt_Hamiltonian-zwei}
for integer exponent $z$ when we put $p=1-e^{-\b}$.
Conversely, summing over $\s$ we get 
\begin{equation*}
\begin{split}
\sum_{\s} K(\sigma,\o)&=C \sum_{\s}\prod_{D \in \D}\Bigl((1-p)1_{\o(D)=0}+ p1_{\o(D)=1}1_D(\s)\Bigr)\cr 
&=C\prod_{D \in \D}(1-p)^{1-\o(D)}p^{\o(D)}q^{k(\o)} 
\end{split}
\end{equation*}
where $k(\o)$ is the number of connected components (in the sense that open $z$-subsets are called connected if they share at least one vertex) of the configuration $\o\in\{0,1\}^\D$ also counting isolated elements of $\D$. 
We call this measure the \textit{generalized random cluster measure} (generalized RCM) assigning probability to configurations of $z$-cliques. More details for the case $z=2$ can be found for example in \cite{GHM00}.

The case $q=1$ is independent percolation on $z$-clique variables since we declare each $z$-clique (subset of $z$ elements) independently to be open with probability $p$ and closed with probability $1-p$. For $q>1$ configurations additionally get $q$-dependent weights which give bias to configurations with many connected components. 

The coupling measure \eqref{Coupling} describes an intimate relation between the generalized Potts measure and the generalized RCM. For example let $C_1,\dots,C_k$ be a partition of $\{1,\dots,N\}$ given by the connected components of a configuration distributed according to the generalized mean-field RCM with parameters $q,z$ and $p=1-e^{-\frac{\b (z-1)!}{N^{z-1}}}$. Then the empirical distribution under the generalized Potts measure with parameters $\b,z$ and $q$ is given by
\begin{equation*}
L_N=\frac{1}{N}\sum_{i=1}^k\a_i|C_i|
\end{equation*}
where the $\a_i$ are independent and equidistributed random variables on $\{\d_1,\dots,\d_q\}$ and we suppressed the additional term in the Hamiltonian \eqref{Alt_Hamiltonian}. Now let us consider the variance of the empirical distribution w.r.t the generalized Potts measure
\begin{equation*}
\begin{split}
\text{Var}_{\pi^N_{\b,q,z}}&[L_N(1)]=\E_{\pi^N_{\b,q,z}}[(L_N(1)-\frac{1}{q})^2]=\E_{\text{RCM}}[(\sum_{i=1}^k\a_i(1)\frac{|C_i|}{N}-\frac{1}{q})^2]\cr
&=\E_{\text{RCM}}[\sum_{j,i=1}^k(\a_i(1)-\frac{1}{q})(\a_j(1)-\frac{1}{q})\frac{|C_i||C_j|}{N^2}]
=\frac{q-1}{q^2}\E_{\text{RCM}}[\sum_{i=1}^k(\frac{|C_i|}{N})^2].\cr
\end{split}
\end{equation*} 
We have $\E_{\text{RCM}}[\max_{i\in\{1,\dots,q\}}(\frac{|C_i|}{N})^2]\leq\E_{\text{RCM}}[\sum_{i=1}^k(\frac{|C_i|}{N})^2]\leq \E_{\text{RCM}}[\max_{i\in\{1,\dots,q\}}(\frac{|C_i|}{N})]$ and hence $\text{Var}_{\pi^N_{\b,q,z}}[L_N(1)]\to0$ iff $\max_{i\in\{1,\dots,q\}}(\frac{|C_i|}{N})\to0$ in probability w.r.t the RCM. In other words, phase-transition of the generalized Potts model is equivalent to percolation of the generalized RCM.

The case $z=2$ has been studied in great detail in \cite{BGJ96}. Under the right scaling $p=\l/N$ the critical value $\l_c$ for percolation of the RCM equals the critical inverse temperature $\b_c$ for phase-transition of the Potts model. We expect the same to be true for the generalized RCM and the generalized Potts measure (on a computational level even for $q$ non-integer valued) with $p=\l/N^{z-1}$.

\medskip
Notice that for the generalized RCM, the assumption of $q$ to be integer-valued can be abandoned. In \cite{B10} again for the case $z=2$ an interesting extension of the Potts measure (on the lattice) the so-called \textit{divide and color model} (DCM) is considered. The DCM is a probability measure on $\{1,\dots,s\}^{\Z^d}$ corresponding to the following two-step procedure: First pick a random edge configuration $\o$ according to the $q$-biased RCM. Secondly assign spin $i\in\{1,\dots,s\}$ independently to every connected component of $\o$ with probability $a_i$ where $\sum_{i=1}^sa_i=1$. For integers $1<s<q$ and $a_i=k_i/q$ with $k_i\in\N$ and $\sum_{i=1}^s k_i=q$ the fuzzy Potts model is contained as a special case. The main result is, that with the exception of the Potts model ($q=s$, $a_i\equiv1/q$) the DCM is Gibbs only for large $p$. Notice that our result about loss of Gibbsianness of the fuzzy Potts model in the low temperature regime is again contained.

\end{document}